\documentclass[a4paper,10pt]{article}
\usepackage{amsmath,amssymb,amsthm}
\usepackage{picins} 
\newtheorem{theorem}{Theorem}[section]

\newtheorem{lemma}[theorem]{Lemma}

\theoremstyle{definition}
\newtheorem{definition}[theorem]{Definition}
\newtheorem{remark}{Remark}

\newcommand{\ep}{\varepsilon}
\newcommand{\eps}[1]{{#1}_{\varepsilon}}
\numberwithin{equation}{section}
\newcommand{\bn}{{\bold n}}
\newcommand{\bz}{{\bold z}}

\newcommand{\bq}{{\bold q}}

\newcommand{\bu}{{\bold u}}
\newcommand{\bv}{{\bold v}}
\newcommand{\bw}{{\bold w}}
\newcommand{\bff}{{\bold f}}
\newcommand{\bg}{{\bold g}}
\newcommand{\bh}{{\bold h}}
\newcommand{\bA}{{\bold A}}

\newcommand{\bD}{{\bold D}}

\newcommand{\bF}{{\bold F}}
\newcommand{\bG}{{\bold G}}

\newcommand{\bI}{{\bold I}}

\newcommand{\bK}{{\bold K}}

\newcommand{\bM}{{\bold M}}

\newcommand{\bQ}{{\bold Q}}
\newcommand{\bS}{{\bold S}}
\newcommand{\bT}{{\bold T}}

\newcommand{\bV}{{\bold V}}

\newcommand{\DV}{{\rm Div}\,}
\newcommand{\dv}{{\rm div}\,}

\newcommand{\BR}{{\Bbb R}}
\newcommand{\BC}{{\Bbb C}}

\newcommand{\BN}{{\Bbb N}}

\newcommand{\hu}{\hat{u}}

\newcommand{\hpi}{\hat{\pi}}
\newcommand{\htheta}{\hat{\theta}}

\newcommand{\CA}{{\mathcal A}}

\newcommand{\CD}{{\mathcal D}}

\newcommand{\CF}{{\mathcal F}}

\newcommand{\CL}{{\mathcal L}}

\newcommand{\CR}{{\mathcal R}}
\newcommand{\CS}{{\mathcal S}}
\newcommand{\CT}{{\mathcal T}}
\newcommand{\CH}{{\mathcal H}}

\newcommand{\CP}{{\mathcal P}}

\newcommand{\CX}{{\mathcal X}}

\newcommand{\pd}{\partial}

\newcommand{\bng}{\bn_{\Gamma(t)}}
\newcommand{\bngz}{\bn_{\Gamma_0}}
\newcommand{\Hol}{{\rm Hol}\,}
\newenvironment{cases*}%
{%
\left\{
\begin{array}{@{}r@{\;}l@{\quad}l@{}}
}%
{\end{array}\right.}

\newcommand{\hbu}{\hat{\bold u}}
\newcommand{\hbw}{\hat{\bold w}}

\newcommand{\hd}{\hat{d}}

\newcommand{\hrho}{\hat{\rho}}

\newcommand{\hmu}{\hat{\mu}}
\newcommand{\hlambda}{\hat{\lambda}}
\newcommand{\hkappa}{\hat{\kappa}}

\begin{document}
\title{Local well-posedness of compressible-incompressible two-phase flows with
phase transitions}
\author{Yoshihiro Shibata
\thanks{Department of Mathematics and Research Institute of 
Science and Engineering \endgraf
Waseda University, 
Ohkubo 3-4-1, Shinjuku-ku, Tokyo 169-8555, Japan. \endgraf
e-mail address: yshibata@waseda.jp \endgraf
Partially supported by JST CREST, JSPS
Grant-in-aid for Scientific Research (S) \# 24224004
and SGU. }}

\date{}
\maketitle
\begin{abstract}
This paper is concerned with the basic model for 
compressible and incompressible two phase flows with phase transitions
The flows are separated by  nearly flat interface 
represented as a graph over the $N-1$ dimensional 
Euclidean space $\BR^{N-1}$ ($N \geq 2$). The local well-posedness is proved by
the Banach fixed point theorem based on the maximal $L_p$-$L_q$ regularity
theorm for the linearized problem. 
\end{abstract}
\section{Introduction}

The three states of matter are 
solids, liquids and gases.  The flow consisting  
of two phases, which are mixed and interacting each other, 
is called the two phase flow. To analyze the two phase 
flow is the important problem in the field of the fluid machine.
For example, nowadays it is known that cavitation noise and 
the damage of hard material for tarbo-machines and ship propellers
are induced by impulsive pressures that are caused by 
the collapse of cloud of bubbles in the water.  In fact, K.~Yamamoto 
\cite{Yamamoto} investigated that the water jets injected from 
submerged nozzle with narrow orifice were observed by a high-speed 
video camera and he  found that 
the several times rebound of the cloud of bubbles creates   
very strong pressure pulses which cause the cavitation noise and 
the damage of hard material. Moreover,  D.~Rosseinelli
et al \cite{Ross} showed the similation of 
cloud cavitation collapse by the high performance computer.
Thus, the study of the cavitation, which is described
by the compressible and incompressible fluid flow mathematically,
 has new depelopment
in the experimental fluid mechanics and computational fluid 
mechanics, rather recently.  On the other hand, the mathematical approach
to two phase problem with liquids and bubbles is rare, 
even when they are described by
the compressible and incompressible viscous fluid flow with sharp interface. 
The author knows only results due to Denisova \cite{Denisova} and 
Kubo, Shibata and Soga \cite{KSS}. In this paper, 
we start with the modeling of the two phase problem which 
can be found in the usual engineering text books without any 
mathematical proof and we prove 
the local well-posedness in the phase transition case,
which has not been yet treated in any mathematical 
literature as far as the author knows.

\section{Modeling.}
Following the Pr\"uss idea  in \cite{PSSS}, we discuss 
the modeling. 
Let $\Omega$ be a domain in the $N$ dimensional Euclidean space
$\BR^N$ ($N\geq2$) with boundary $\Gamma_0$.
 Let $\Omega_-$ be a subdomain
of $\Omega$ with boundary $\Gamma$. We assume that 
$\Gamma = \pd\Omega_- \subset \Omega$ and that $\Gamma_0
\cap \Gamma =\emptyset$. Set $\Omega_+ = 
\Omega-\overline{\Omega_-}$. Let  $\varphi=\varphi(\xi, t)
= (\varphi_1(\xi, t), \ldots, \varphi_N(\xi, t))$ be a function 
defined on the closure of $\Omega$ for each time variable $t \in (0, T)$,
$\xi = (\xi_1, \ldots, \xi_N)$ being the reference coordinate system.  
We assume that the map $\xi \to \varphi(\xi, t)$ is one to
one for each $t \in (0, T)$
\footnote[2]{Since this subsection is concerned with
the modeling,
we do not care the regularity of boundary and the map 
$\varphi$. Moreover, we do not mention any integrability of 
functions regorously.  These are formulated mathematically 
in sections 2.}. 
Set $(\pd_t\varphi)(\xi, t) = \bv(x, t)$ with $x = \varphi(\xi, t)$, 
\begin{align*}
\Omega_\pm(t) &= \{x = \varphi(\xi, t) \mid \xi \in \Omega_\pm\}, \quad
\Gamma(t)  = \{x = \varphi(\xi, t) \mid \xi \in \Gamma\}, 
\end{align*}
and $\dot\Omega(t) = \Omega_-(t) \cup \Omega_+(t)$. 
Let $\bng$ be the unit outer normal to $\Gamma(t)$ pointed
from $\Omega_-(t)$ to $\Omega_+(t)$ and let 
$\bngz$ the unit outer normal to $\Gamma_0$. 
Set 
$$[[v]] =v_- - v_+ \quad\text{(the jump of $v$ accross $\Gamma(t)$)}$$
for any $v$ defined on $\dot\Omega(t)$.  Here and hereafter, 
we write $v_\pm = v|_{\Omega_\pm(t)}$. Moreover, given $v_\pm$ 
defined on $\Omega_\pm(t)$, we define $v$ by
$v(x) = v_\pm(x)$ for $x \in \Omega_\pm(t)$. Let 
$H_\Gamma = -\dv_\Gamma\bn_\Gamma$ be the mean curvature of 
$\Gamma(t)$. 
During our modeling, we use the well-known 
Reynolds transport theorem:
$$\frac{d}{dt} \int_{\dot\Omega(t)}f\,dx
= \int_{\dot\Omega(t)}\pd_t f\,dx + \int_{\Gamma(t)}
[[f]]\bv\cdot\bng\,d\nu
+ \int_{\Gamma_0}f\bv\cdot\bngz\,d\nu,$$ 
where $d\nu$ represents the surface element not only of
$\Gamma(t)$ but also of $\Gamma_0$. 
In this orientation, we know that
\begin{equation}\label{0}
\frac{d}{dt}|\Gamma(t)| = -\int_{\Gamma(t)}H_\Gamma
\bv\cdot\bn_\Gamma\,d\nu
\end{equation}
Here and in the following, 
we use  bold small letters to denote $N$-vector and $N$-vector valued
functions and the bold capital letters to denote $N\times N$ matrix and 
$N\times N$ matrix valued functions, respectively. For 
$\bv=(v_1, \ldots, v_N)$ and $\bw = 
(w_1, \ldots, w_N)$, we set $<\bv, \bw> = \bv\cdot\bw = \sum_{j=1}^N
v_jw_j$, which is the usual inner product of $\BR^N$.

In the following, we use the following notation:
\begin{itemize}
\item $\rho: \dot\Omega(t) \to \BR_+$ 
is the mass field, 
\item $\bu: \dot\Omega(t) \to \BR^N$ the velocity field,
\item $\pi : \dot\Omega(t) \to \BR$ the pressure field, 
\item $\bT : \dot\Omega(t) \to \{ \bA \in GL_N(\BR) \mid {}^T\bA = \bA\}$
the stress tensor field,
\item $\theta : \dot\Omega(t) \to \BR_+$ the thermal
field,
\item $e: \dot\Omega(t) \to \BR$
the internal energy, 
\item $\bq: \dot\Omega(t) \to \BR^N$ the heat flux, 
\item $\eta: \dot\Omega(t) \to \BR$ the entropy, 
\end{itemize}
where,  we have set $\BR_+ = (0, \infty)$.  
For the modeling, we use the following Navier-Stokes-Fourier 
system of equations: for $x \in \dot\Omega(t)$ ad $t > 0$ 
\begin{align}
&\pd_t\rho + \dv(\rho\bu) = 0;  \label{1}\\
&\pd_t(\rho\bu) + \dv(\rho\bu\otimes\bu) - 
\dv\bT = 0;\label{2}\\
&\pd_t(\frac{\rho}{2}|\bu|^2 + \rho e)
+ \dv((\frac{\rho}{2}|\bu|^2 + \rho e)\bu)
-\dv(\bT\bu - \bq) = 0.\label{3}
\end{align}
Here, for any $\bu = (u_1, \ldots, u_N)$,
$\bu\otimes\bu$ is the $N\times N$ matrix
whose $(i, j)$ component is  $u_iu_j$, and
for any $\bw = (w_1, \ldots, w_N)$
and  $\bS = (S_{ij})$ 
the divergence forms $\dv\bw$ and $\dv\bS$ are
defined by
$$\dv \bw = \sum_{j=1}^N\pd_jw_j, \enskip
\dv \bS = (\sum_{j=1}^N \pd_j S_{1j}, \ldots, 
\sum_{j=1}^N \pd_j S_{Nj}).$$
During our discussion of the jump condition on 
$\Gamma(t)$ and boundary condition on $\Gamma_0$. 
we assume that $\bv \not=\bu$ on $\Gamma(t)$,
but $\bv = \bu$ on $\Gamma_0$. 

First, we consider the mass conservation:
\begin{equation}\label{mass:1}
\frac{d}{dt}\int_{\dot\Omega(t)} \rho\,dx=0.
\end{equation}
 By \eqref{1} 
 and the Reynolds transport theorem, we have 
\begin{align*}
\frac{d}{dt}\int_{\dot\Omega(t)} \rho\,dx &
 = \int_{\dot\Omega(t)} \pd_t \rho\,dx
+ \int_{\Gamma(t)} [[\rho]]\bv\cdot\bng\,d\nu
+ \int_{\Gamma_0} \rho\bu\cdot\bngz\,d\nu\\
& = -\int_{\dot\Omega(t)} \dv(\rho\bu)\,dx 
+ \int_{\Gamma(t)} [[\rho]]\bv\cdot\bng\,d\nu
+ \int_{\Gamma_0} \rho\bu\cdot\bngz\,d\nu \\
& =-\int_{\Gamma(t)} [[\rho(\bu - \bv)]]\cdot\bng\,d\nu.
\end{align*}
Thus, to obtain \eqref{mass:1}, it is sufficient to
assume that 
\begin{equation}\label{4}
[[\rho(\bu - \bv)]]\cdot\bng = 0 \quad
\text{on $\Gamma(t)$}. 
\end{equation}
In this case, $\rho_+(\bu_+-\bv)\cdot\bng
= \rho_-(\bu_--\bv)\cdot\bng$ on $\Gamma(t)$, so that
the phase flux $\j$ is defined by 
\begin{equation}\label{pf} \j = \rho_+(\bu_+-\bv)\cdot\bng
= \rho_-(\bu_--\bv)\cdot\bng.
\end{equation} 
\begin{itemize}
\item When $\j=0$,  we have  
$[[\bu]]\cdot\bn_\Gamma = 0$. 
\item When $\j\not=0$ and $[[\rho]]\not=0$, we have
\begin{equation}\label{5}
\j = \frac{[[\bu]]\cdot\bng}{[[1/\rho]]}.
\end{equation}
\item When $\j\not=0$ and $[[\rho]]=0$, $\j$ can not 
be decided by the velocity field $\bu$.
\end{itemize}
The case where $\j=0$ is called 
without phase transition and the case where $\j\not=0$  
with phase transition. 


Next, we consider the conservation of momentum:
\begin{equation}\label{momentum:1}
\frac{d}{dt}\int_{\dot\Omega(t)} \rho\bu\,dx =0.
\end{equation}
By \eqref{2} and the Reynolds transport theorem,
we have  
\begin{align*}
&\frac{d}{dt}\int_{\dot\Omega(t)} \rho\bu\,dx =
-\int_{\dot\Omega(t)} \dv(\rho\bu\otimes\bu)\,dx + 
\int_{\dot\Omega(t)} \dv\bT\,dx \\
& = 
-\int_{\Gamma(t)} ([[\rho\bu\otimes(\bu-\bv) -[[\bT]])
\bng\,d\nu + \int_{\Gamma_0} 
\bT\bngz\,d\nu.
\end{align*}
Thus, in order that $\frac{d}{dt}\int_{\dot\Omega(t)} \rho\bu\,dx=0$
holds, it is sufficient to assume that 
\begin{equation}\label{6}\begin{cases*}
[[\rho\bu\otimes(\bu-\bv) - \bT]]\bng &= \dv_\Gamma\bT_\Gamma
&\quad\text{on $\Gamma(t)$}, \\
\bT\bngz &= 0
&\quad\text{on $\Gamma_0$}
\end{cases*}\end{equation}
Here, $\bT_\Gamma$ is the stress tensor field
on $\Gamma(t)$. Note that $\int_{\Gamma(t)}\dv_\Gamma T_\Gamma
\,d\nu =0$. 
We assume  that $\dv_\Gamma\bT_\Gamma = -\sigma H_\Gamma\bng$, 
where $\sigma$ is a  non-negative constant
 describing the coefficient of 
surface tension.

We represent the interface condition \eqref{6} with the help of 
the phase flux $\j$ as follows: 
\begin{align*}
[[\rho\bu\otimes(\bu-\bv)]]\cdot\bng = \rho_-\bu_-(\bu_--\bv)\cdot\bng
-\rho_+\bu_+(\bu_+-\bv)\cdot\bng
= \j[[\bu]].
\end{align*}
Moreover, by \eqref{1} we rewrite \eqref{2} as follows: 
$$\pd_t(\rho\bu) + \dv(\rho\bu\otimes\bu)
= \bu(\pd_t\rho + \dv(\rho\bu)) + \rho(\pd_t\bu
+ \bu\cdot\nabla\bu)
= \rho(\pd_t\bu
+ \bu\cdot\nabla\bu).$$
Summing up, we have obtained 
\begin{equation}\label{6*}
\begin{cases*}
\rho(\pd_t\bu + \bu\cdot\nabla\bu) -
\dv\bT &= 0 &\quad\text{in $\dot\Omega(t)$}, \\
\j[[\bu]]-[[\bT\bng]] &= -\sigma H_\Gamma\bng
&\quad\text{on $\Gamma(t)$}, \\
\bT\bngz &= 0
&\quad\text{on $\Gamma_0$}.
\end{cases*}\end{equation}
Here and in the following, for any $N$-vector functions 
$\bw=(w_1, \ldots, w_N)$, $\bz = (z_1, \ldots, z_N)$
and scalor function $f$, we set 
$\bw\cdot\nabla f = \sum_{j=1}^N w_j\pd_jf$, and 
$\bw\cdot\nabla\bz$ is an $N$-vector function whose $i$ th component is
$\bw\cdot\nabla z_i$. 


Next, we consider the balance of energy. 
We look for a sufficient condition to obtain  
the conservation of energy: 
\begin{equation}\label{energy:1}
\frac{d}{dt}(\int_{\dot\Omega(t)}(\frac{\rho}{2}|\bu|^2 + \rho e)\,dx
+ \sigma|\Gamma(t)|)= 0.  
\end{equation}
By \eqref{3} and the Reynolds transport theorem, we have  
\begin{align*}
&\frac{d}{dt}\int_{\dot\Omega(t)}(\frac{\rho}{2}|\bu|^2 + \rho e)\,dx \\
& = -\int_{\Gamma(t)} [[(\frac{\rho}{2}|\bu|^2 + \rho e)(\bu-\bv)
-(\bT\bu - \bq)]]\cdot\bng\,d\nu 
+ \int_{\Gamma_0} (\bT\bu - \bq)\cdot\bngz\,d\nu,
\end{align*}
which, combined with \eqref{0}, furnishes that  
\begin{align*}
&\frac{d}{dt}(\int_{\dot\Omega(t)}(\frac{\rho}{2}|\bu|^2 + \rho e)\,dx
+ \sigma|\Gamma(t)|) 
= \int_{\Gamma_0} (\bT\bu - \bq)\cdot\bngz\,d\nu \\
&\quad
-\int_{\Gamma(t)}([[(\frac{\rho}{2}|\bu|^2 + \rho e)(\bu-\bv)
-(\bT\bu - \bq)]]\cdot\bng + \sigma H_\Gamma \bv\cdot\bng)\,d\nu.  
\end{align*}
Thus, in order to obtain \eqref{energy:1}, 
 it is sufficient to assume that 
\begin{equation}\label{8}
\begin{cases*}
[[(\frac{\rho}{2}|\bu|^2 + \rho e)(\bu-\bv)
-(\bT\bu - \bq)]]\cdot\bng + \sigma H_\Gamma \bv\cdot\bng&=0
&\quad\text{on $\Gamma(t)$}, \\
(\bT\bu - \bq)\cdot\bngz &= 0 
&\quad\text{on $\Gamma_0$}.
\end{cases*}\end{equation}
Since $\bT\bngz=0$ on $\Gamma_0$, we assume that 
$\bq\cdot\bngz=0$ on $\Gamma_0$. 
By \eqref{pf} and \eqref{6*}, 
\begin{align*}
&[[\frac{\rho}{2}|\bu|^2(\bu-\bv)]]\cdot\bng \\
&\quad = \frac{\j}{2}(|\bu_--\bv+\bv|^2 - |\bu_+-\bv+\bv|^2) 
= \frac{\j}{2}[[|\bu-\bv|^2]] + \j[[\bu]]\cdot\bv\\
& \quad = \frac{\j}{2}[[|\bu-\bv|^2]] + [[\bT\bng]]\cdot\bv
-\sigma H_\Gamma \bv\cdot\bng\end{align*}
Since $[[\rho e(\bu-\bv)]]\cdot\bng = \j[[e]]$, the first
equation of \eqref{8} becomes: 
$$
\frac{\j}{2}[[|\bu-\bv|^2]] +\j[[e]] 
- [[\bT(\bu-\bv)]]\cdot\bng 
+ [[\bq]]\cdot\bng =0.
$$
Moreover, using \eqref{1} and \eqref{2}, we rewrite \eqref{3}
as follows: 
\begin{align*}
&\pd_t(\frac{\rho}{2}|\bu|^2 + \rho e)+
\dv((\frac{\rho}{2}|\bu|^2 + \rho e)\bu) - \dv
(\bT\bu - \bq) \\
&= (\frac12|\bu|^2+e)\pd_t\rho + \rho(\bu\cdot\pd_t\bu
+ \pd_t e) + (\frac12|\bu|^2 + e)\dv(\rho\bu)
+ \rho\bu\cdot(\bu\cdot\nabla\bu + \nabla e)\\
&\quad- (\dv\bT)\cdot \bu -\bT:\nabla\bu + \dv\bq
\\
&= (\frac12|\bu|^2+e)(\pd_t\rho + \dv(\rho\bu))
+ \bu\cdot(\rho (\pd_t\bu + \bu\cdot\nabla\bu) - \dv\bT)\\
&\quad+ \rho(\pd_t e+ \bu\cdot\nabla\bu) - \bT:\nabla\bu + \dv\bq.
\end{align*}
Here, we have set $\bT:\nabla\bu = \sum_{i,j=1}^NT_{ij}\pd_iu_j$.
Thus, we have 
$$\rho(\pd_t e + \bu\cdot\nabla e) + \dv\bq - 
\bT:\nabla\bu = 0.$$
Summing up, we have obtained 
\begin{equation}\label{9}\begin{cases*}
\rho(\pd_t e + \bu\cdot\nabla e) + \dv\bq - 
\bT:\nabla\bu &= 0 &\quad\text{in $\dot\Omega(t)$}, \\
\frac{\j}{2}[[|\bu-\bv|^2]] +\j[[e]] 
- [[\bT(\bu-\bv)]]\cdot\bng 
+ [[\bq]]\cdot\bng&=0&\quad\text{on $\Gamma(t)$}, \\
\bq\cdot\bngz&= 0 &\quad\text{on $\Gamma_0$}.
\end{cases*}\end{equation}

The number of interface conditions is so far  not enough. To find one 
more condition, we consider  the law of entropy increase:
\begin{equation}\label{entropy:1}
\frac{d}{dt} \int_{\dot\Omega(t)} \rho\eta\,dx\geq 0.
\end{equation}
 For this purpose,
we introduce the constitutive  laws in the phases. 
According to the Newton law,  the  stress tensor $\bT$ is given by 
$$\bT = 2\mu\bD(\bu) + (\lambda-\mu)\dv\bu\bI- \pi\bI.$$
Here, $\bD(\bu) = \frac12({}^T\nabla \bu + \nabla\bu)$ is 
the strain tensor field, $\bI$ the $N\times N$ identity matrix,
$\mu$ and  $\lambda$ are the first and second 
viscosity coefficients satisfying the condition:
\begin{equation}\label{11} \mu > 0, \quad 
\lambda \geq \frac{N-2}{N}\mu.
\end{equation}
To prove local well-posedness, it suffices to 
assume that $\mu > 0$ and $\lambda > 0$.
According  to the Fourier law, the heat flux $\bq$ is given by 
\begin{equation}\label{12}
\bq = -d\nabla\theta.
\end{equation}
with thermal conductivity $d$ satisfying the condition: $d > 0$. 
Moreover, the first law of thermodynamics tells us that the internal energy 
$e$, the entropy $\eta$, and 
the pressure term $\pi$  have the  relation: 
\begin{equation}\label{13}
de = \theta d\eta + \frac{\pi}{\rho^2}d\rho.
\end{equation}
We define the free energy $\psi$ for the unit mass  and 
the specific heat $\kappa$ by 
\begin{equation}\label{14}
\psi = e - \theta\eta, \quad \kappa
 = \frac{\pd e}{\pd\theta},
\end{equation}
respectively.  We assume that $\kappa > 0$. 
Since $\dfrac{\pd e}{\pd \eta} = \theta$ and $\dfrac{\pd e}{\pd \rho}
= \dfrac{\pi}{\rho^2}$ as follows from \eqref{13}, by \eqref{1} 
\begin{equation}\label{eneg:0}\begin{split}
\pd_te + \bu\cdot\nabla e &= \frac{\pd e}{\pd \eta}(\pd_t \eta
+ \bu\cdot\nabla\eta) + \frac{\pd e}{\pd \rho}
(\pd_t \rho+ \bu\cdot\nabla\rho) \\
& = \theta(\pd_t \eta
+ \bu\cdot\nabla\eta)
- \frac{\pi}{\rho}\dv\bu.
\end{split}\end{equation}
In addition, 
\begin{equation}\label{eneg:1}
\bT:\nabla\bu = 2\mu|\bD(\bu)|^2 + (\lambda-\mu)(\dv\bu)^2 - \pi\dv\bu,
\end{equation}
which, combined with the first equation of \eqref{9}, 
\eqref{eneg:0}, and \eqref{eneg:1}, furnishes
\begin{equation}\label{ent:1}
\rho\theta(\pd_t\eta + \bu\cdot\nabla\eta) - \dv(d\nabla\theta)
-(2\mu|\bD(\bu)|^2+(\lambda-\mu)(\dv\bu)^2) = 0.
\end{equation}
On the other hand, we have 
\begin{equation}\label{ent:2}
\pd_t(\rho\eta) + \dv(\rho\eta\bu)
= \eta(\pd_t\rho+\dv(\rho\bu))+\rho(\pd_t\eta+\bu\cdot\nabla\eta)
= \rho(\pd_t\eta + \bu\cdot\nabla\eta).
\end{equation}
In the following, we assume that 
\begin{equation}\label{heat:1}
\theta > 0, \quad [[\theta]]=0.
\end{equation}
Since $\theta$ represents the absolute temperature, $\theta>0$ is 
natural assumption. While phase transition happens, 
the temperature does not change, so that $[[\theta]]=0$ is also natural 
assumption. By \eqref{ent:1} and \eqref{ent:2}
\begin{equation}\label{21}
\pd_t(\rho\eta) + \dv(\rho\eta\bu)
= \frac{1}{\theta}\{\dv(d\nabla\theta) 
+2\mu|\bD(\bu)|^2 + (\lambda-\mu)
(\dv\bu)^2\}.
\end{equation}
By the Reynolds transport theorem, \eqref{21},  
and the divergence theorem of Gauss,  
\begin{align*}
&\frac{d}{dt}\int_{\dot\Omega(t)}\rho\eta\,dx\\
& = -\int_{\dot\Omega(t)}\dv(\rho\eta\bu)\,dx
+ \int_{\dot\Omega(t)}\frac{1}{\theta}
\{\dv(d\nabla\theta)
+(2\mu|\bD(\bu)|^2 + (\lambda-\mu)(\dv\bu)^2\}\,dx
\\
&\quad +\int_{\Gamma(t)}[[\rho\eta]]\bv\cdot\bng\,d\nu
+\int_{\Gamma_0}\rho\eta\bu\cdot\bngz\,d\nu\\
& = -\int_{\Gamma(t)}[[(\rho\eta)(\bu-\bv)]]\cdot\bng\,d\nu
+ \int_{\Gamma(t)}[[\frac{d}{\theta}\nabla\theta]]\cdot\bng\,d\nu
+ \int_{\Gamma_0} \frac{1}{\theta}(d\nabla\theta\cdot\bngz)\,d\nu
\\
& \quad + \int_{\dot\Omega(t)}\{\frac{d|\nabla\theta|^2}{\theta^2}
+ 2\mu|\bD(\bu)|^2 + (\lambda-\mu)(\dv\bu)^2\}\,dx.
\end{align*}
Since $2\mu|\bD(\bu)|^2 + (\lambda-\mu)(\dv\bu)^2 \geq 0$
as follows from \eqref{11}, to obtain \eqref{entropy:1}
 it is  sufficient to assume that 
\begin{equation}\label{step:1}\begin{split}
[[(\rho\eta)(\bu-\bv)-\frac{d}{\theta}\nabla\theta]]\cdot\bng =0
&\quad\text{on $\Gamma(t)$}, \\ 
d\nabla\theta\cdot\bngz = 0 &\quad\text{on $\Gamma_0$}.
\end{split}\end{equation}
Moreover, by $[[\theta]]=0$ and \eqref{pf}, 
 the first equation of \eqref{step:1} becomes 
\begin{equation}\label{23}
\j[[\theta\eta]] - [[d\nabla\theta]]\cdot\bng = 0
\quad\text{on $\Gamma(t)$}, 
\end{equation}
which is called the Stefan law.
 In fact, by \eqref{pf} and $[[\theta]]=0$, 
\begin{align*}
0 & = [[(\rho\eta)(\bu-\bv)-\frac{d}{\theta}\nabla\theta]]\cdot\bng\\
&= (\rho_-\eta_-)(\bu_--\bv)\cdot\bng - (\rho_+\eta_+)(\bu_+-\bv)\cdot\bng
- \frac{(d_-\nabla\theta_--d_+\nabla\theta_+)\cdot\bng}{\theta}\\
&= \frac{1}{\theta}(\j[[\theta\eta]] - [[d\nabla\theta]]\cdot\bng).
\end{align*}
Note that the Stefan law becomes the usual 
jump condition: $[[d\nabla\theta]]\cdot\bng = 0$ on $\Gamma(t)$
provided that $\j=0$.

Next, assuming that $\j\not=0$ and $[[\rho]]\not=0$ and 
using \eqref{23}, we rewrite \eqref{9}. Given $\bw$, we set 
$\CT_{\bng}\bw = \bw - (\bw\cdot\bng)\bng$, which is
the tangential part of $\bw$ along $\bng$. Since  
$\bw = (\bw\cdot\bng)\bng + \CT_{\bng}\bw$ and 
$\CT_{\bng}\bw \cdot \bng = 0$, we have  
\begin{equation}\label{orth:1}
|\bw|^2 = |\bw\cdot\bng|^2 + |\CT_{\bng}\bw|^2.
\end{equation} 
In the following, we assume that
\begin{equation}\label{10}
\CT_{\bng}[[\bu - \bv]]=0.
\end{equation}
Especially, we have 
\begin{equation}\label{10*}
[[\bu]]= ([[\bu]]\cdot\bng)\bng.
\end{equation}
Since $|\bw|^2 - |\bv|^2 = (\bw-\bv)\cdot(\bw+\bv)$ for any 
$\bw$ and $\bv$, by \eqref{10} $[[|\CT_{\bng}(\bu-\bv)|^2]] = 0$,
so that by \eqref{pf} and \eqref{orth:1}   
\begin{equation}\label{gibbs:1}
[[\frac12|\bu-\bv|^2]] = \frac12[[|(\bu-\bv)\cdot\bng|^2]]
= \frac{\j^2}{2}[[\frac{1}{\rho^2}]].
\end{equation}
Since $e = \psi + \theta\eta$, we have $[[e]]=[[\psi]]+[[\theta\eta]]$.
Thus, recalling \eqref{12} and using  \eqref{23} and \eqref{gibbs:1}, 
we rewrite the jump condition in \eqref{9} as follows: 
\begin{align*}
0 &=\frac{\j}{2}[[|\bu-\bv|^2]] + \j[[e]] - [[\bT(\bu-\bv)]]\cdot\bng
 - [[\bq]]\cdot\bng\\
& = \j^3[[\frac{1}{2\rho^2}]] - [[\bT(\bu-\bv)]]\cdot\bng + 
\j[[\psi]].
\end{align*}
Moreover, noting that $\bT_\pm$ are symmetric matrices, 
we have 
\begin{align*}
[[\bT(\bu-\bv)]]\cdot\bng&= [[(\bu-\bv)\cdot\bT\bng]] \\
&= [[((\bu-\bv)\cdot\bng)\bng\cdot\bT\bng]]
+[[\CT_{\bng}(\bu-\bv)\cdot\bT\bng]].
\end{align*}
By \eqref{pf}, 
$$
[[((\bu-\bv)\cdot\bng)\bng\cdot\bT\bng]]
= \j[[\frac{1}{\rho}\bng\cdot\bT\bng]].
$$
On the other hand, by  \eqref{6*} and
\eqref{10} and \eqref{10*}, 
\begin{multline*}
[[\CT_{\bng}(\bu-\bv)\cdot\bT\bng]]
= \CT_{\bng}(\bu_--\bv)\cdot[[\bT\bng]]\\
 = \CT_{\bng}(\bu_--\bv)\cdot(\j[[\bu]]\cdot\bng + \sigma H_\Gamma)
\bng= 0.
\end{multline*}
Thus, we have obtained 
$0 = \j([[\psi]] + \j^2[[\dfrac{1}{2\rho^2}]]
-[[\dfrac{1}{\rho}\bng\cdot\bT\bng]])$. 
Since $\j\not=0$, finally we arrive at the condition: 
\begin{equation}\label{25}
[[\psi]] + \j^2[[\frac{1}{2\rho^2}]]-
[[\frac{1}{\rho}\bng\cdot\bT\bng]]=0 \quad\text{on $\Gamma(t)$},
\end{equation}
which is called the generalized Gibbs-Thomson law. 

Finally, we calculate $V_\Gamma := \bv\cdot\bng$. By \eqref{pf}
we have $\bv\cdot\bng = \bu_-\cdot\bng -\dfrac{\j}{\rho_-}$. 
When $\j=0$, it follows from \eqref{pf} and \eqref{10} that $[[\bu]]=0$, 
so that 
$\bv\cdot\bng = \bu\cdot\bng$.
When $\j\not=0$ and $[[\rho]]\not=0$, by \eqref{5} we have 
$\j = \dfrac{[[\bu]]\cdot\bng}{[[1/\rho]]}$, so that 
\begin{align*}
\bv\cdot\bng
& = \bu_-\cdot\bng -\frac{\j}{\rho_-}
 = \frac{[[\rho\bu]]\cdot\bng}{[[\rho]]}.
\end{align*}

Summing up, we have obtained 
\begin{alignat}2
V_\Gamma : = \bv\cdot\bng &= \bu\cdot\bng &\quad&(\j=0), \nonumber\\
V_\Gamma := \bv\cdot\bng &= \frac{[[\rho\bu]]\cdot\bng}{[[\rho]]}
&\quad&(\j\not=0\text{ and }[[\rho]]\not=0).
\label{26}
\end{alignat}

Next, we consider the case where 
$\j\not=0$ and $[[\rho]]=0$. In this case,  $[[\bu]]\cdot\bng=0$,
which, combined with \eqref{10},  furnishes that 
$[[\bu]]=0$,  so that 
\eqref{6*} becomes
\begin{equation}\label{26*}
[[\bT\bng]] = \sigma H_\Gamma\bng \quad\text{on
$\Gamma(t)$}.
\end{equation}
To derive \eqref{25}, we assume that
$[[\rho]]\not=0$, so that we reconsider the second condition of
\eqref{9}. By 
$[[\bu]]=0$,  $[[|\bu-\bv|^2]]=0$. By \eqref{12} and \eqref{23},
 $\j[[e]] +[[\bq]]\cdot\bng = \j[[\psi]]$. 
In addition, by \eqref{pf}, \eqref{26*}, and the symmetricity of $\bT_\pm$
\begin{multline*}
[[\bT(\bu-\bv)]]\cdot\bng = (\bu_--\bv)\cdot\bT_-\bng - 
(\bu_+-\bv)\cdot\bT_-\bng \\
= (\bu_--\bv)\cdot[[\bT\bng]] 
= (\bu_--\bv)\cdot\sigma H_\Gamma\bng 
 = \j\frac{\sigma}{\rho_-}H_\Gamma. 
\end{multline*}
Since $\j\not=0$, the second equation of \eqref{9} becomes  
\begin{equation}\label{27}
[[\psi]] - \frac{\sigma}{\rho_-}H_\Gamma = 0. 
\quad\text{on $\Gamma(t)$}.
\end{equation}

Noting that $\pd_t\rho + \bu\cdot\nabla\rho = -\rho\dv\bu$ as follows 
from \eqref{1} and recalling the formulas: $\dfrac{\pd e}{\pd\theta} = \kappa$
and $\dfrac{\pd e}{\pd\rho} = \dfrac{\pi}{\rho^2}$
(cf. \eqref{13} and \eqref{14}), 
we have 
\begin{align*}
\rho(\pd_te+\bu\cdot\nabla e)
&= \rho(\frac{\pd e}{\pd\theta}\pd_t\theta
+ \frac{\pd e}{\pd\rho}\pd_t\rho 
+ \frac{\pd e}{\pd\theta}\bu\cdot\nabla\theta
+  \frac{\pd e}{\pd\rho}\bu\cdot\nabla\rho) \\
&= \rho\kappa(\pd_t\theta + \bu\cdot\nabla\theta)
-\frac{\pi}{\rho}\dv\bu,
\end{align*}
which, combined with the first eqution in \eqref{9} and \eqref{eneg:1},
furnishes that 
$$
\rho\kappa(\pd_t\theta + \bu\cdot\nabla\theta)
-\dv(d\nabla\theta) 
-(2\mu|\bD(\bu)|^2+ (\lambda-\mu)(\dv\bu)^2)
+\pi(1-\frac{1}{\rho})\dv\bu = 0.
$$

Summing up, we have obtained the equations: for $x \in \dot\Omega(t)$
and $t > 0$ 
\begin{align}
\pd_t\rho + \dv(\rho\bu) &= 0, \nonumber  \\
\rho(\pd_t\bu + \bu\cdot\nabla\bu)-\dv\bT
&= 0, \label{29}\\
\rho\kappa(\pd_t\theta + \bu\cdot\nabla\theta)-
\dv(d\nabla\theta)&=
(2\mu|\bD(\bu)|^2+ (\lambda-\mu)(\dv\bu)^2)
-\pi(1-\frac{1}{\rho})\dv\bu, \nonumber
\end{align}
subject to the boudary condition: for $x\in\Gamma_0$ 
and $t > 0$:
\begin{equation}\label{30}
\bT\bngz = 0, \quad
d\nabla\theta\cdot\bngz = 0\quad\text{on $\Gamma_0$},
\end{equation}
and one of the following interface conditions: 
for $x \in\Gamma(t)$ and $t > 0$ 
\vskip0.5pc\noindent
\thetag1 When $\j = 0$, 
\begin{equation}\label{31}\begin{split}
&[[\bu]]=0, \quad [[\bT\bng]] = \sigma H_\Gamma\bng,\quad 
[[\theta]] = 0, \\ 
&[[d\nabla\theta\cdot\bng]]=0, \quad
\bv\cdot\bng = \bu\cdot\bng.
\end{split}\end{equation}
\thetag2 When $\j \not= 0$ and $[[\rho]]\not=0$, 
\begin{equation}\label{32}\begin{split}
&\CT_{\bng}[[\bu]]=0, \quad \j[[\bu]]
-[[\bT\bng]] = -\sigma H_\Gamma\bng,\quad [[\theta]] = 0, \\
&\j[[\theta\eta]] -[[d\nabla\theta\cdot\bng]]=0, \quad
[[\psi]]+\j^2[[\frac{1}{2\rho^2}]] - [[\frac{1}{\rho}
\bng\bT\bng]]= 0,\\
&\bv\cdot\bng =\frac{[[\rho\bu]]\cdot\bng}{[[\rho]]}, 
\quad
\j = \frac{[[\bu]]\cdot\bng}{[[1/\rho]]}.
\end{split}\end{equation}
\thetag3 When $\j \not= 0$ and $\rho = \rho_- = \rho_+$ (constants), 
\begin{equation}\label{33}\begin{split}
&[[\bu]]=0, \quad 
[[\bT\bn_\Gamma]] = \sigma H_\Gamma\bng,\quad 
[[\theta]] = 0, \quad
\j[[\theta\eta]] -[[d\nabla\theta\cdot\bng]]=0, \\
&\rho[[\psi]]-\sigma H_\Gamma = 0,
\quad \bv\cdot\bng = \bu\cdot\bng - \frac{\j}{\rho}.
\end{split}\end{equation} 
\begin{remark} Assuming that $\Omega_-= \Omega$ and $\Omega_+=\emptyset$,
we have the one phase problem.  In this case, as boundary condition on
$\Gamma_0$, we have 
$$\bT\bngz = \sigma H_{\Gamma_0}\bng, \quad d\nabla\theta\cdot\bngz = 0 
\quad\text{on $\Gamma_0$}.
$$
\end{remark}
\section{Problem}
The problem of this paper is concerned with 
 the compressible and incompressible two phase flow separated
by a nearly flat interface with  phase transition.
Let $h_0(x')$ be a given function with respect to $x' = (x_1, 
\ldots, x_{N-1})$ and we set 
\begin{align*}
\Omega_\pm &= \{x = (x_1, \ldots, x_N) \in \BR^N
\mid \pm (x_N - h_0(x')) > 0 \quad
\text{for $x' \in \BR^{N-1}$}\}, \\
\Gamma &= \{x \in \BR^N \mid x_N = h(x') \quad 
\text{for $x' \in \BR^{N-1}$}\}.
\end{align*}
In this case, $\Omega = \BR^N$ and $\Gamma_0 = \emptyset$. 
Let $h(x', t)$ be a unknown function and we assume that
the time evolutions of 
domains $\Omega_\pm$ and  the surface $\Gamma$ are given by
\begin{equation}\label{interface:1}\begin{split}
\Omega_\pm(t) &= \{x = (x_1, \ldots, x_N) \in \BR^N
\mid \pm (x_N - h(x', t)) > 0 \quad 
\text{for $x' \in \BR^{N-1}$}\}, \\
\Gamma(t) &= \{x \in \BR^N \mid x_N = h(x',t) \quad
\text{for $x' \in \BR^{N-1}$}\}.
\end{split}\end{equation}
In this case,  $\bng = (-\nabla'h, 1)/\sqrt{1+|\nabla'h|^2}$,
with $\nabla'h=(\pd_1h, \ldots, \pd_{N-1}h)$ ($\pd_j = \pd/\pd x_j$). 
Moreover, $\varphi(x, t) = (x', x_N+h(x', t))$, so that 
$\bv\cdot\bng = \pd_t\varphi\cdot\bng = \pd_th/\sqrt{1 + |\nabla'h|^2}$.

In view of \eqref{eneg:1},
 \eqref{29}, \eqref{30} and \eqref{32}, our problem 
is given as follows: 
\begin{equation}\label{eq:1}\begin{split}
&\text{For $x \in \Omega_+(t)$, $t>0$},\\
&\begin{cases}
\rho_+(\pd_t\bu_++ \bu_+\cdot\nabla\bu_+) - \DV\bT_+  = 0, \quad 
\pd_t\rho_+ + \dv(\rho_+\bu_+)  = 0, \\[0.5pc]
\rho_+ \kappa_+
(\pd_t\theta_+ + \bu_+\cdot\nabla\theta_+) 
-\dv(d_+\nabla\theta_+)- \bT_+:\nabla\bu_+
-\dfrac{\pi}{\rho}\dv\bu_+ = 0
\end{cases}
\\
&\text{and, for $x \in \Omega_-(t)$, $t>0$},\\
&\begin{cases}
&\rho_{*-}(\pd_t\bu_-+ \bu_-\cdot\nabla\bu_-) - \DV\bT_-  = 0, \quad 
\dv\bu_-  = 0, \\
&\rho_{*-} \kappa_-
(\pd_t\theta_- + \bu_-\cdot\nabla\theta_-) 
-\dv(d_-\nabla\theta_-)- \bT_-:\nabla\bu_- = 0
\end{cases}
\end{split}
\end{equation}
subject to the jump conditions: for $x \in \Gamma(t)$ and $t > 0$,
\begin{equation}\label{jump:1}
\begin{cases}
[[\dfrac1\rho]]\j^2\bn_\Gamma - [[\bT\bn_\Gamma]] 
= -\sigma H_\Gamma\bn_\Gamma, 
&\quad \CT_{\bng}[[\bu]]  =0, \\[0.5pc]
\j[[\theta\eta]] - [[d\dfrac{\pd\theta}{\pd\bn_\Gamma}]]  = 0, &\quad
[[\theta]] = 0, \\[0.5pc]
[[\psi]] + [[\dfrac{1}{2\rho^2}]]\j^2 - [[\frac{1}{\rho}
\bn_\Gamma\cdot\bT\bn_\Gamma]]  = 0, &\quad
\pd_t h = \dfrac{[[\rho\bu]]\cdot(-\nabla'h, 1)}{[[\rho]]},\\[0.8pc]
\quad
\j=\dfrac{[[\rho\bu]]\cdot\bng}{[[\rho]]},
\end{cases}
\end{equation}
and the initial conditions:
\begin{equation}\label{initial}\begin{split}
&(\rho_+, \bu_+, \theta_+)|_{t=0} = (\rho_{*+}+\rho_{0+}, \bu_{0+},
\theta_{*} + \theta_{0+}) \enskip\text{in $\Omega_+$}, \\
&(\bu_-, \theta_-)|_{t=0} = (\bu_{0+},
\theta_{*} + \theta_{0+}) \enskip\text{in $\Omega_-$}, 
\quad 
h|_{t=0} = h_0 \enskip\text{on $\Gamma$}.
\end{split}\end{equation}
Here, $\rho_{*\pm}$, $\theta_{*}$ and $\sigma$ 
 are positive constants describing the 
reference mass densities of $\Omega_\pm$, the reference temperature
of $\Omega_\pm$ and the coefficient of 
the surface tension, respectively. Moreover, 
 $\bT_\pm = \bS_\pm - \pi_\pm\bI$ with
\begin{align*}
&\bS_+ =\bS_+(\bu_+, \rho_+, \theta_+))
= \mu_+\bD(\bu_+) + (\lambda_+-\mu_+)\dv\bu\bI, \\
&\bS_-=\bS_-(\bu_-, \theta_-) = \mu_-\bD(\bu_-).
\end{align*}
Here, $d_+ =d_+(\rho, \theta)$, $\mu_+ = \mu_+(\rho, \theta)$, 
$\lambda_+ = \lambda_+(\rho, \theta)$, 
$\kappa_+ = \kappa_+(\rho, \theta)$
 are positive 
$C^\infty$ functions with respect to $(\rho, \theta) 
\in (0, \infty)\times(0, \infty)$, and 
$\psi_+(\theta, \rho)$ and $\eta_+(\theta, \rho)$ 
are  real valued $C^\infty$
functions with respect to $(\rho, \theta) \in (0, \infty)\times(0, \infty)$,
while  $d_- =d_-(\theta)$, $\mu_- = \mu_-(\theta)$, 
$\kappa_- = \kappa_-(\theta)$
 are positive 
$C^\infty$ functions with respect to $\theta 
\in (0, \infty)$, and 
$\psi_-(\theta)$ and $\eta_-(\theta)$ are real valued $C^\infty$
functions with respect to $\theta \in (0, \infty)$.
Moreover, 
we also assume that $\pi_+$ is given by $\pi_+ = P_+(\rho, \theta)$,
where $P_+$ is some $C^\infty$ function
with respect to $(\rho, \theta) \in (0, \infty)\times(0, \infty)$
such that $\dfrac{\pd P_+}{\pd\rho} > 0$ for any 
$(\rho, \theta) \in (0, \infty)\times(0, \infty)$. 

The main purpose of this paper is to show  
the local wellposedness of 
 problem \eqref{eq:1}, \eqref{jump:1} and \eqref{initial}
in the maximal $L_p$-$L_q$ regularity class
under the assumption that $\rho_{*\pm}$ and 
$\theta_*$ satisfy the condition:
\begin{equation}\label{assump:0}
\rho_{*-} \not= \rho_{*+}, \quad
\psi_-(\theta_*) - \psi_+(\rho_{*+}, \theta_*)
+\Bigl(\frac{1}{\rho_{*-}}-\frac{1}{\rho_{*+}}\Bigr)
P_+(\rho_{*+}, \theta_*) = 0.
\end{equation}
To state our main result, we transform $\Gamma(t)$ to the 
flat interface. 
Set 
$$
\BR^N_\pm = 
\{x=(x_1, \ldots, x_N) \in \BR^N \mid \pm x_N > 0\}, \quad
\BR^N_0 = \{x \in \BR^N \mid x_N = 0\}.
$$
We transfer the problem given in domains $\Omega_\pm(t)$
to that in $\dot\BR^N = \BR^N_+\cup \BR^N_-$ with 
interface $\BR^N_0$. Let $h(x',t)$ be 
a function appearing in the definition of $\Gamma(t)$ 
in \eqref{interface:1}.
Let $H(x, t)$ be a solution 
to the equations:
$(1-\Delta)H = 0$ in $\BR^N_\pm$ with $H|_{x_N=0} = h(x', t)$, where 
 $\Delta H = \sum_{j=1}^N \pd_j^2H$. To prove the local well-posedness,
we assume that $h_0$ is small enough, so that we may assume that
\begin{equation}\label{assump:2}
1 + \frac{\pd }{\pd x_N}H(x,t) \geq \frac12 
\quad\text{for any $x \in \BR^N_\pm$
and $t \in (0, T)$}.
\end{equation}
If we consider the transformation: 
\begin{equation}\label{trans:1}
y_N = x_N + H(x, t), \quad 
y_j = x_j \enskip (j = 1, \ldots, N-1), 
\end{equation} 
then by \eqref{assump:2} $\Omega_\pm(t)$ and $\Gamma(t)$ 
are transformed to $\BR^N_\pm$ and 
$\BR^N_0$, respectively, 
because $y_N = h(y',t)$ when $x_N=0$ and $\frac{\pd y_N}{\pd x_N}
= 1 + (\frac{\pd H}{\pd x_N})(x, t) \geq \frac12$.

  Let $\bu_\pm$, $\rho_+$,
$\pi_-$ 
and $\theta_\pm$ satisfy problem \eqref{eq:1}, \eqref{jump:1} and 
\eqref{initial}. Set
\begin{align*}
\hat\bu_\pm(x, t) & = \bu_\pm(x', x_N+H(x, t), t), \quad
\hat\rho_+(x, t)  = \rho_+(x', x_N+H(x, t), t),\\
\hat\pi_-(x, t) & = \pi_-(x', x_N+H(x, t), t)-
P_+(\rho_{*+}, \theta_*), \\
\hat\theta_\pm(x, t)&= \theta_\pm(x', x_N+H(x, t), t),\quad
\mu_{*+} = \mu_+(\rho_{*+}, \theta_{*}), \quad 
\lambda_{*+} = \lambda_+(\rho_{*+}, \theta_{*}), \\
\kappa_{*+} & = \kappa_+(\rho_{*+}, \theta_{*}), \quad
d_{*+} = d_+(\rho_{*+}, \theta_{*}), \quad 
\mu_{*-} = \mu_-(\theta_{*}), \quad 
\kappa_{*-} = \kappa_-(\theta_{*}), \\
\quad d_{*-} &= d_-(\theta_{*}),\quad
\hmu_+  = \mu_+(\hrho_+, \htheta_+), \quad 
\hlambda_+ = \lambda_+(\hrho_+, \htheta_+), \quad
\hmu_- = \mu_-(\htheta_-),\\
\hkappa_+ &= \kappa_+(\hrho_+, \htheta_+), \quad
\hkappa_- = \hkappa_-(\htheta_-), \quad
\hat d_+ = d_+(\hrho_+, \theta_+), 
\quad \hat d_- = d_-(\htheta_-),\\
\tilde\rho_+ &= \hrho_+ - \rho_{*+}, \quad  
\tilde\mu_\pm  = \hmu_\pm - \mu_{*\pm}, \quad 
\tilde\lambda_+ = \hlambda_+ - \lambda_{*+}, \quad 
\tilde\kappa_\pm  = \hkappa_\pm - \kappa_{*\pm}, \\
\tilde d_\pm &= \hat d_\pm - d_{*\pm}.
\end{align*}
Setting 
$H_0 = \pd_tH$, $H_j = \pd_jH$ ($j=1, \ldots, N$), we have 
\begin{equation}\label{change:1}\begin{split}
(\pd_tf)(x', x_N+H(x, t), t) &= \pd_t\hat f(x, t) - \frac{H_0}{1+H_N}
\pd_N\hat f(x, t), \\
(\pd_jf)(x', x_N+H(x, t), t) &= \pd_j\hat f(x, t) - \frac{H_j}{1+H_N}
\pd_N\hat f(x, t) \quad(j=1, \ldots, N).
\end{split}\end{equation}
In the following, we set
$$K_j = \dfrac{H_j}{1+H_N} \enskip(j = 0, 1, \ldots, N), \quad 
\bK = (K_1, \ldots, K_N),
\quad \bK_0 = (K_0, \bK).$$
By \eqref{change:1} we have 
$$
\nabla \pi_- = \bQ\nabla\hpi_- = \left(\begin{matrix}
1 & 0& \cdots &0 & K_1 \\
\vdots & & &\vdots & \vdots \\
0&0&\cdots&1&K_{N-1} \\
0&0&\cdots &0 & \frac{1}{1+H_N}
\end{matrix}\right)
\left(\begin{matrix}
\pd_1\hpi_-\\\vdots\\  \\
 \pd_N\hpi_-
\end{matrix}\right), 
$$
and $\bQ^{-1}$ is given by  
$$
\bQ^{-1} = \left(\begin{matrix}
1 & 0 & \cdots & 0 & - H_1 \\
\vdots & & &\vdots & \vdots \\
0&0&\cdots & 1 & - H_{N-1} \\
0&0&\cdots &0& 1 + H_N
\end{matrix}\right)
= \bI + \bQ_1 \quad\text{with}\enskip
\bQ_1 = \left(\begin{matrix}
0&\cdots&0&-H_1 \\
\vdots &\ddots & \vdots & \vdots \\
0&\cdots&0& -H_{N-1} \\
0&\cdots&0& H_N
\end{matrix}\right).
$$
By \eqref{change:1} we  have 
\begin{equation}\label{div:1}\begin{split}
\dv\bu_\pm &= \dv \hbu_\pm + V_{\dv}(\hbu_\pm, H) \\
&=\frac{1}{1+H_N}\bigl\{\dv\hbu_\pm - f_-(\hbu_\pm, H)\bigr\} \\
&= \frac{1}{1+H_N}\bigl\{\dv\hbu_\pm - \dv\bff_-(\hbu_\pm, H)\bigr\}
\end{split}\end{equation}
with 
\begin{align*}
&V_\dv(\hbu_\pm, H) = -\sum_{j=1}^NK_j\pd_N\hu_{\pm j}, \quad 
f_-(\hbu_\pm, H) = \sum_{j=1}^{N-1}
(H_N\pd_j\hu_{\pm j} - H_j\pd_N\hu_{\pm j}), \quad\\
&\bff_-(\hbu_\pm, H) = -(H_N\hu_{\pm 1}, \ldots, 
H_N\hu_{\pm N-1}, -\sum_{j=1}^{N-1}H_j\hu_{\pm j}).
\end{align*}
For any $N\times N$ matrix of functions $\bG = 
{}^T(\bg_1, \ldots, \bg_N)$, where ${}^T\bM$ denotes the 
transposed $\bM$, with ${}^T\bg_i = (g_{i1}, \ldots, g_{iN})$,
by \eqref{div:1} we have
\begin{equation}\label{div:1*}
\DV \bG = \DV\hat\bG + \bV_{\DV}(\hat\bG, H)
\end{equation}
with $\bV_{\DV}(\hat\bG, H) = 
{}^T(V_\dv(\hat\bg_1, H), \ldots, 
V_\dv(\hat\bg_N, H))$. 
Moreover, we set 
\begin{equation}\label{div:2}
D_{ij}(\bu_\pm) = D_{ij}(\hbu_\pm) + V_{D_{ij}}(\hbu_\pm, H), 
\enskip \bD(\bu_\pm) = \bD(\hbu_\pm) + \bV_\bD(\hbu_\pm, H),
\end{equation}
where 
$V_{D_{ij}}(\hbu_\pm, H) = -(K_i\pd_N\hu_{\pm j} + 
K_j\pd_N\hu_{\pm i})$ and 
$\bV_{\bD}(\hbu_\pm, H)$ is the $N\times N$ matrix whose 
$(i, j)$ component is $V_{D_{ij}}(\hbu_\pm, H)$.

Under these preparations, we see easily that problem
\eqref{eq:1}, \eqref{jump:1} and \eqref{initial} is transformed
to the following problem:
\begin{alignat}2
&\left\{
\begin{aligned}
&\pd_t\hrho_+ + \bv_+\cdot\nabla\hrho_+ + 
\hrho_+(\dv\hbu_+ + V_\dv(\hbu_+, H)) = 0 \\
&\rho_{*+}\pd_t\hbu_+ - \DV\bS_{*+}(\hbu_+) = \bF_+\\
&\rho_{*+}\kappa_{*+}\pd_t\htheta_+ - d_{*+}\Delta\htheta_+ 
= F_{\theta +}
\end{aligned}\right.
&\quad&\text{in $\BR^N_+\times(0, T)$, } \nonumber \\
&\left\{
\begin{aligned}
&\rho_{*-}\pd_t\hbu_- - \DV\bS_{*-}(\hbu_-) + \nabla\hpi_- 
= \bF_-\\
&\dv\hbu_- = f_- = \dv\bff_-\\
&\rho_{*-}\kappa_{*-}\pd_t\htheta_- - d_{*-}\Delta\htheta_- 
= F_{\theta -}
\end{aligned}\right.
&\quad&\text{in $\BR^N_-\times(0, T)$,}\nonumber\\
&\left\{
\begin{aligned}
&\mu_{*-}D_{iN}(\hbu_-)|_- - \mu_{*+}D_{iN}(\hbu_+)|_+ = G_i \\
&T_- - T_+ =\sigma \Delta'H+ G_N, \quad \rho_{*-}^{-1}T_-
-\rho_{*+}^{-1}T_+ = G_{N+1} \\
&\hu_{-i}|_- - \hu_{+i}|_+ = K_i \\
&\htheta_{-}|_- - \htheta_{+}|_+ =0,  \quad
d_{*-}\pd_N\htheta_-|_- - d_{*+}\pd_N\htheta|_+ = 
G_\theta\\
&\pd_tH - \Bigl(\dfrac{\rho_{*-}}{\rho_{*-}-\rho_{*+}}\hu_{-N}|_-
-\dfrac{\rho_{*+}}{\rho_{*-}-\rho_{*+}}\hu_{+N}|_+
\Bigr) = G_h\end{aligned}\right. 
&\quad&\text{on $\BR^N_0\times(0, T)$,} \nonumber \\
&\left\{
\begin{aligned}
&(\hrho_+, \hbu_+, \htheta_+)|_{t=0}
 = (\hrho_{0+}, 
\hbu_{0+}, \htheta_{0+}) \quad \text{in $\BR^N_+$},
 \\
&(\hbu_-, \htheta_-)|_{t=0}
 = (\hbu_{0-}, \htheta_{0-}) \quad \text{in $\BR^N_-$},
\quad
 H|_{t=0} = H_0 \quad\text{on $\BR^N_0$}
\end{aligned}\right.
\label{stokeseq:1}
\end{alignat}
where $i=1, \ldots, N-1$, and we have set 
\begin{align}
&T_- = 
(\mu_{*-}D_{NN}(\hbu_-)-\hpi_-)|_-, \\
&T_+ = 
(\mu_{*+}D_{NN}(\hbu_+)
+ (\lambda_{*+}-\mu_{*+})\dv\hbu_+), \nonumber \\
&\bv_+ = (\hu_{+1}, \ldots, \hu_{+N-1}, 
\hu_{+N}-K_0-\sum_{j=1}^NK_j\hu_{+j}), \nonumber \\
& \bS_{*+}(\bu) = \mu_{*+}\bD(\bu) + 
(\lambda_{*+} - \mu_{*+})\dv\bu\bI,\quad
\bS_{*-}(\bu) = \mu_{*-}\bD(\bu), \nonumber \\
&\hrho_{0+}(x) = \rho_{0+}(x', x_N+H_0(x)), \enskip
\hbu_{0\pm}(x) = \bu_{0\pm}(x', x_N+H_0(x)), \nonumber \\
&\htheta_{0\pm}(x) = \theta_{0\pm}(x', x_N+H_0(x)).
\label{quant:1}
\end{align}
Here, $f|_\pm(x_0) = \lim_{x\to x_0 \atop x \in \Omega_\pm(t)}
f(x)$ for $x_0 \in \BR^N_0$.  Moreover, $H_0$ is a function 
satisfying the equations $(1-\Delta)H_0 = 0$ in 
$\BR^N_\pm$ and $H_0|_{x_N=0} = h_0$, and 
the right-hand sides in \eqref{stokeseq:1}
 are defined by the following 
formulas: 
\allowdisplaybreaks{ 
\begin{align*}
&\bF_+=\bF_+(\hrho_+, \hbu_+, H) =\\
&-\tilde\rho_+\pd_t\hbu_+ 
+\hrho_+\{ K_0\pd_N\hbu_+ 
- \hbu_+\cdot\nabla\hbu_+ 
+ (\hbu_+\cdot\bK)\pd_N\hbu_+\} \\
&+ \DV(\tilde\mu_+\bD(\hbu_+) + \hmu_+\bV_\bD(\hbu_+, H))
+\bV_\bD(\hmu_+(\bD(\hbu_+) + \bV_\bD(\hbu_+, H)))\\
&+\nabla\{(\tilde\lambda_+ - \tilde\mu_+)\dv\hbu_+ +
(\hlambda_+ - \hmu_+)V_{\dv}(\hbu_+, H)\}\\ 
&+\bV_{\DV}((\hlambda_+-\hmu_+)
(\dv\hbu_+ + V_{\dv}(\hbu_+, H))\bI) 
-\bQ\nabla P_+(\hrho_+, \htheta_+), \\
&F_{\theta +} = F_{\theta +}(\hrho_+, \hbu_+, \htheta_+, H) 
= -d_{*+}\sum_{j=1}^N\pd_j(K_j\pd_N\htheta_+)\\
&+\sum_{j=1}^N\pd_j(\tilde d_+(\pd_j\htheta_+
 -K_j\pd_N\htheta_+)) 
+ \sum_{j=1}^NK_j\pd_N(\hd_+(\pd_j\htheta_+ - K_j
\pd_N\htheta_+))\\ 
&
-(\hrho_+\hkappa_+-\rho_{*+}\kappa_{*+})
\pd_t\htheta_+ 
+ \hrho_+\hkappa_+(K_0\pd_N\htheta_+ 
- \hbu_+\cdot\nabla\htheta_+
+ (\hbu_+\cdot\bK)\pd_N\htheta_+) \\
&
 + 2\hmu_+|\bD(\hbu_+) + \bV_\bD(\hbu_+, H)|^2
 + (\hlambda_+-\hmu_+)
(\dv\hbu_+ + V_\dv(\hbu_+, H))^2 \\
& +P_+(\hrho_+, \htheta_+) 
(1-\frac{1}{\hrho_+})
(\dv\hbu_++V_\dv(\hbu_+, H)), \\
&\bF_-= \bF_-(\hbu_-, H) =
-\bQ_1(\rho_{*-}\pd_t\hbu_- - \mu_{*-}\DV\bD(\hbu_-))\\
&-(I+\bQ_1)\{\rho_{*-}(K_0\pd_N\hbu_- - \hbu_-\cdot\nabla\hbu_- +
(\hbu_-\cdot\bK)\pd_N\hbu_-) \\
&
+\DV(\tilde\mu_-\bD(\hbu_-) + \hmu_-\bV_\bD(\hbu_-, H))
+ \bV_\DV(\hmu_-(\DV\bD(\hbu_-) + 
\bV_\bD(\hbu_-, H)))\},\\ 
&f_=f_-(\hbu_-, H) = \sum_{j=1}^{N-1}
\{H_N\pd_j\hu_{-j} - H_j\pd_N\hu_{-j}), \\
&\bff_- = \bff_-(\hbu_-, H) = 
-(H_N\hu_{-1}, \ldots, H_N\hu_{-\,N-1}, -\sum_{j=1}^{N-1}
H_j\hu_{-j}), \\
&F_{\theta -} = F_{\theta -}(\hbu_-, \htheta_-, H) 
=-\rho_{*-}\tilde\kappa_-\pd_t\htheta_-
-d_{*-}\sum_{j=1}^N\pd_j(K_j\pd_N\htheta_-)\\
&+ \rho_{*-}\hkappa_-(K_0\pd_N\htheta_- 
- \hbu_-\cdot\nabla\htheta_-
+ (\hbu_-\cdot\bK)\pd_N\htheta_-) 
+ \sum_{j=1}^N\pd_j(\tilde d_+(\pd_j\htheta_-
-K_j\pd_N\htheta_-)) \\
&
+ \sum_{j=1}^NK_j\pd_N(\hd_-(\pd_j\htheta_- - K_j
\pd_N\htheta_-))
 + 2\hmu_-|\bD(\hbu_-) + \bV_\bD(\hbu_-, H)|^2,\\
&G_i = G_i(\hrho_+, \hbu_\pm, H) =\\
& -\{(\tilde\mu_-D_{iN}(\hbu_-) + \hmu_-V_{D_{iN}}(\hbu_-, H))|_-
- (\tilde\mu_+D_{iN}(\hbu_+) + \hmu_+V_{D_{iN}}(\hbu_+, H))|_+\} \\
&+\sum_{j=1}^{N-1}(\pd_jH)\{\hmu_-(D_{ij}(\hbu_-)
+V_{D_{ij}}(\hbu_-, H))|_- - 
\hmu_+(D_{ij}(\hbu_+)
+V_{D_{ij}}(\hbu_+, H))|_+\}\\
&
-(\pd_iH)[\sum_{j=1}^{N-1}(\pd_jH)\{\hmu_-
(D_{ij}(\hbu_-) + V_{D_{ij}}(\hbu_-))|_-
- \hmu_+
(D_{ij}(\hbu_+) + V_{D_{ij}}(\hbu_+))|_+\}, \\
&-\{\hmu_-(D_{NN}(\hbu_-) + V_{D_{NN}}(\hbu_-, H))|_-
-\hmu_+(D_{Nj}(\hbu_+) + V_{D_{Nj}}(\hbu_+, H))|_+\}],\\
&G_N = G_N(\hrho_+, \hbu_\pm, H)=\\
&-(\tilde\mu_-D_{NN}(\hbu_-) + \hat\mu_-V_{D_{NN}}(\hbu_-, H))|_-
-(\tilde\mu_+D_{NN}(\hbu_+) + \hat\mu_+V_{D_{NN}}(\hbu_+, H))|_+\\
&+\{(\tilde\lambda_+-\tilde\mu_+)\dv\bu_+ + (\hlambda_+
-\hmu_+)V_{\dv}(\hbu_+, H) -(P_+(\hrho_+, \htheta_+)-
P_+(\rho_{*+}, \theta_{*+}))\}|_+\\
&+\sum_{j=1}^{N-1}(\pd_jH)(
\hat\mu_-(D_{Nj}(\hbu_-) + V_{D_{Nj}}(\hbu_-, H))|_- 
- \hat\mu_+(D_{Nj}(\hbu_+) + V_{D_{Nj}}(\hbu_+, H))|_+)\\
&+\sigma\Bigl\{(1-\frac{1}{\sqrt{1+|\nabla'H|^2}}\Delta'H
-\sum_{i,j=1}^{N-1}\frac{(\pd_iH)(\pd_i\pd_jH)}
{(1 + |\nabla'H|^2)^{3/2}}\Bigr\}\\
&+\Bigl(\frac{1}{\rho_{*-}}-\frac{1}{\hat \rho_+|_+}\Bigr)^{-1}
(\hat u_{-N}|_- - \hat u_{+N}|_+)^2
(1 + |\nabla'H|^2), \\
&G_{N+1} = G_{N+1}(\hrho_+, \hbu_\pm, \htheta_\pm, H)=\\
&-\frac{1}{\rho_{*-}}(P_+(\hrho_+, \htheta_+)-P_+(\rho_{*+},
\theta_{*+}))-\Bigl(\frac{1}{\hrho_+|_+}-\frac{1}{\rho_{*+}}\Bigr)
P_+(\hrho_+, \htheta_+) \\
&+(\psi_-(\htheta_-) - \psi_-(\theta_*))|_-
-(\psi_+(\hrho_+, \htheta_+) - \psi_+(\rho_{*+}, \theta_*))|_+
\\
&+(\psi_-(\htheta_-)|_- - \psi_+(\hrho_+, \htheta_+)|_+)
|\nabla'H|^2 \\
&+\frac12\Bigl(\frac{1}{\rho_{*-}}+\frac{1}{\hat\rho_+|_+}
\Bigr)\Bigl(\frac{1}{\rho_{*-}}-\frac{1}{\hat\rho_+|_+}
\Bigr)^{-1}(u_{-N}|_- -u_{+N}|_+)^2(1 + |\nabla'H|^2)^2\\
&-\Bigl\{\frac{1}{\rho_{*-}}\tilde\mu_-D_{NN}(\hbu_-)|_-
-\frac{1}{\rho_{*+}}(\tilde\mu_+D_{NN}(\hbu_+)
+ (\tilde\lambda_+-\tilde\mu_+)\dv\bu_+)|_+\Bigr\}\\
&+\Bigl(\frac{1}{\hrho_+|_+}-\frac{1}{\rho_{*+}}\Bigr)
\{\hmu_+(D_{NN}(\hbu_+) + V_{D_{NN}}(\hbu_+, H)) \\
&\phantom{+(\frac{1}{\hrho_+|_+}-\frac{1}{\rho_{*+}})
\{\hmu_+(D_{NN}(\hbu_+)+\,}
+(\hlambda_+-\hmu_+)(\dv\hbu_+ + V_{\dv}(\hbu_+, H))\}|_+\\
&-\sum_{i,j=1}^{N-1}\Bigl\{\frac{\hmu_-}{\rho_{*-}}
(D_{ij}(\hbu_-) + V_{D_{ij}}(\hbu_-, H))|_-\\
&\phantom{+(\frac{1}{\hrho_+|_+}-\frac{1}{\rho_{*+}})
\{\hmu_+(D_{NN}}
-\frac{\hmu_+}{\hrho_+}(D_{ij}(\hbu_+)+V_{D_{ij}}(\hbu_+, H))|_+
\Bigr\}(\pd_iH)(\pd_jH)\\
&+2\sum_{i=1}^{N-1}\Bigl\{
\frac{\hmu_-}{\rho_{*-}}(D_{iN}(\hbu_-) + 
V_{D_{iN}}(\hbu_-,H))|_-\\
&\phantom{+(\frac{1}{\hrho_+|_+}-\frac{1}{\rho_{*+}})
\{\hmu_+(D_{NN}}
-\frac{\hmu_+}{\hrho_+}(D_{iN}(\hbu_+)+V_{D_{iN}}(\hbu_+, H))|_+
\Bigr\}(\pd_iH)\\
&-|\nabla'H|^2\Bigl[\sigma
\Bigl(\frac{\Delta'H}{\sqrt{1 + |\nabla'H|^2}}
-\sum_{i,j=1}^{N-1}\frac{(\pd_iH)(\pd_i\pd_j H)}
{(1 + |\nabla'H|^2)^{3/2}}\Bigr)\\
&+\Bigl(\frac{1}{\rho_{*-}}-\frac{1}{\hrho_+|_+}\Bigr)^{-1}
(u_{-N}|_- - u_{+N}|_+)^2(1+|\nabla'H|^2)\\
& +\sum_{i=1}^{N-1}\{\hmu_-(D_{iN}(\hbu_-)
+V_{D_{iN}}(\hbu_-, H))|_- 
-\hmu_+(D_{iN}(\hbu_+)+ V_{D_{iN}}(\hbu_+, H))|_+
\}(\pd_iH)\\
&-\{\hmu_-D_{NN}(\hbu_-) + V_{D_{NN}}(\hbu_+, H))|_-
-\hmu_+(D_{NN}(\hbu_+) + 
V_{D_{NN}}(\hbu_+, H))|_+\},\\
&K_i = K_i(\hbu_\pm, H) = -(\pd_iH)(\hu_{-N}|_- - \hu_{+N}|_+), \\
&G_\theta = G_\theta(\hrho_+, \hbu_\pm, \htheta_\pm, H)= \\
&(1 + |\nabla'H|^2)(\hu_{-N}|_- - \hu_{+N}|_+)
\Bigl(\frac{1}{\rho_{*-}}
-\frac{1}{\hrho_+|_+}\Bigr)^{-1}
(\htheta_-\eta_-(\htheta_-)|_-
-\htheta_+\eta_+
(\hrho_+, \htheta_+)|_+)\\
&-(\tilde d_-(\nabla\htheta_- - \bK\pd_N\htheta_-)|_-
-\tilde d_+(\nabla\htheta_+ - \bK\pd_N\htheta_+)|_+)\cdot(-\nabla'H, 1)\\
&+ (d_{*-}\nabla'\htheta_-|_- - d_{*+}\nabla'\htheta_+|_+)
\cdot\nabla'H
+(d_{*-}\pd_N\htheta_-|_- - d_{*+}\pd_N\htheta_+|_+)
\bK\cdot(-\nabla'H, 1), \\
&G_h = G_h(\hrho_+, \hbu_\pm, H) =\\
&\Bigl(\frac{1}{\rho_{*-}-\hrho_+|_+}-
\frac{1}{\rho_{*-}-\rho_{*+}}\Bigr)
(\rho_{*-}\hat u_{-N}|_-- \rho_{*+}\hat u_{+N}|_+)
+ \frac{\rho_{*+}-\hrho_+}{\rho_{*-}-\hrho_+}
\hat u_{+N}|_+ \\
& - \frac{\rho_{-*}}{\rho_{*-}-\hrho_+|_+}\sum_{j=1}^{N-1}
(\pd_jH)\hat u_{-j}|_-
+ \frac{\hrho_+}{\rho_{*-}-\hrho_+}\sum_{j=1}^{N-1}
(\pd_jH)\hat u_{+j}|_+.
\end{align*}
}
The phase flux $\j$ is eliminated by using the formula:
$$\j = (\hu_{-N}|_- - \hu_{+N}|_+)\Bigl(\frac{1}{\rho_{*-}}
-\frac{1}{\hrho_+|_++\rho_{*+}}\Bigr)^{-1}
\sqrt{1 + |\nabla'H|^2} \quad\text{on $\BR^N_0\times(0, T)$}.$$
Moreover, we have used the formulas: for $x \in \Gamma(t)$ and 
$t > 0$
\begin{align*}
&\pi_--P_++(\lambda_+-\mu_+)\dv\bu_+ \\
&=-\sigma H_\Gamma -[[\frac1\rho]]\j^2
-\sum_{j=1}^{N-1}[[\mu D_{Nj}]](\pd_jH) 
+ [[\mu D_{NN}]],  \\
&H_\Gamma\bng = \Bigl\{
\dv'\Bigl(\frac{\nabla'H}
{1 + |\nabla'H|^2}\Bigr)\Bigr\}(-\nabla'H, 1)
/\sqrt{1 + |\nabla'H|^2}
\end{align*}
where $\nabla'H = (\pd_1H, \ldots, \pd_{N-1}H)$ and $\dv' \bv' = 
\sum_{j=1}^{N-1}\pd_jv_j$ for $\bv' = (v_1, \ldots, v_{N-1})$.

Since $T_--T_+ = \sigma\Delta'H + G_N$ and 
$\rho_{*-}^{-1}T_- - \rho_{*+}^{-1}T_+ = G_{N+1}$ are equivalent to
\begin{equation}\label{eqiv:1}
T_\pm = \frac{\rho_\pm\sigma}{\rho_{*-} - \rho_{*+}}\Delta'H
+ \frac{\rho_{*-}\rho_{*+}}{\rho_{*-}-\rho_{*+}}
(\rho_{*\mp}^{-1}G_N - G_{N+1}),
\end{equation}
the compatibility condition for problem \eqref{stokeseq:1} is 
\begin{align}
&\dv\hbu_{-0} = f_-(\hbu_{-0}, H_0) = \dv\bff_-(\hbu_{-0}, H_0)
\quad\text{in $\BR^N_-$}, \nonumber\\
&\mu{*-}D_{iN}(\hbu_{-0})|_- - \mu_{*+}D_{iN}(\hbu_{+0})|_+ = 
G_i(\hrho_{0+}, \hbu_{0\pm}, H_0)
\quad(i=1, \ldots, N-1), \nonumber \\
&\hu_{0-i}|_- - \hu_{0+i}|_+ = K_i(\hbu_\pm, H_0) 
\quad(i=1, \ldots, N-1), \nonumber \\
&\htheta_{0-}|_- - \htheta_{0+}|_+ =0, \nonumber \\
&d_{*-}\pd_N\htheta_{0-}|_- - d_{*+}\pd_N\htheta_{0+}|_+ = 
G_\theta(\hrho_{0+}, \hbu_{0\pm}, \htheta_{0\pm}, H_0), 
\nonumber \\
&(\mu_{*+}D_{NN}(\hbu_{0+})
+ (\lambda_{*+}-\mu_{*+})\dv\hbu_{0+})|_+ = 
\frac{\rho_+\sigma}{\rho_{*-} - \rho_{*+}}\Delta'h_0
\nonumber \\
&\quad
+ \frac{\rho_{*-}\rho_{*+}}{\rho_{*-}-\rho_{*+}}
(\rho_{*-}^{-1}G_N(\hrho_+, \hbu_{0\pm}, H_0)
 - G_{N+1}(\hrho_+, \hbu_{0\pm}, \htheta_{0\pm} ,H_0).
\label{compat:1}
\end{align}
The following theorem is the main result of this
paper concerning  
the local well-posedness of problem \eqref{stokeseq:1}.
\begin{theorem}\label{main:2}
Let $1 < p, q < \infty$ with $2/p + N/q < 1$. 
Assume that $\rho_{*\pm}$ and $\theta_{*}$ satisfy the condition
\eqref{assump:0}. 
Then, given any positive time $T$, there exists an $\epsilon > 0$ such 
that problem \eqref{stokeseq:1} admits
unique  solutions $\hrho_+$, $\hbu_\pm$ and $\htheta_\pm$
with 
\begin{align*}
&\hrho_+ \in W^1_p((0, T), L_q(\BR^N_+)) \cap L_p((0, T),W^1_q(\BR^N_+)),
\\
&(\hbu_\pm, \htheta_\pm) \in W^1_p((0, T), L_q(\BR^N_\pm))
\cap L_p((0, T), W^2_q(\BR^N_\pm)),\\
& H \in W^1_p((0, T), W^2_q(\BR^N)) \cap L_p((0, T), W^3_q(\BR^N))
\end{align*}
for any initial data 
$$\hrho_{0+} \in W^1_q(\BR^N_+), \enskip  
(\hbu_{0\pm}, \htheta_{0\pm}) 
\in B^{2(1-1/p)}_{q,p}(\BR^N_\pm), \enskip H_0 \in W^{3-1/p}_{q,p}(\BR^N)$$
satisfying  the smallness condition:
$$\|\hrho_{0+}\|_{W^1_q(\BR^N_+)} + \sum_{\ell=\pm}
\|(\hbu_{0\ell}, \htheta_{0\ell})\|_{B^{2(1-1/p)}_{q,p}(\BR^N_\ell)}
+ \|H_0\|_{B^{3-1/p}_{q,p}(\BR^N)} \leq \epsilon
$$
and compatibility condition \eqref{compat:1}.
\end{theorem}
Here and in the following, $L_q(\BR^N_\pm)$ and $W^m_q(\BR^N_\pm)$ 
denote the usual Lebesgue space and Sobolev space of order $m$
in the $L_q(\BR^N_\pm)$ sense, while $\|\cdot\|_{L_q(\BR^N_\pm)}$
and $\|\cdot\|_{W^m_q(\BR^N_\pm)}$ denote their norms, repsectively.
For the Banach space $X$, $L_p((0, T), X)$ and $W^m_p((0, T), X)$ 
denote the Lebesgue space and the Sobolev space with values in 
$X$, while $\|\cdot\|_{L_p((0, T), X)}$ and $\|\cdot\|_{W^m_p((0, T), X)}$
denote their norms, respectively.  $B^{(1-\theta)a + \theta b}_{q, p}
(\BR^N_\pm) $ denotes the real interpolation space $(W^a_q(\BR^N_\pm), 
W^b_q(\BR^N_\pm))_{\theta, p}$ with real interpolation functor 
$(\cdot, \cdot)_{\theta, p}$ and $0 < \theta < 1$.
\begin{remark}
\thetag1 The mathematical study of the compressible and 
incompressible two phase problem is quite rare as far as 
the author knows.  First Denisova \cite{Denisova} studied the evolution 
of the compressible and incompressible two phase flow 
with sharp interface without phase transition under
some restriction on the viscosity coefficients.
Recently, Kubo, Shibata and Soga \cite{KSS} studied the same problem as 
in \cite{Denisova} without any restriction on viscosity 
coefficients in case of without surface tension and 
without phase transition.  This paper is 
the first manuscript to treat the compressible and 
incompressible two phase problem with phase transition
\footnote{Modeling and the main results of this paper were 
announced in the abstract of $39^{\rm th}$ Sapporo symposium
on PDE at Hokkaido University (cf. \cite{S3}).}. 
The incompressible and incompressible two phase problem
with phase transition was studied by J.~Pr\"uss, 
et al. \cite{PSSS, PS, PSW}. 
\end{remark}

\section{Maximal $L_p$-$L_q$ regularity} \label{sec:3}

In the following, we assume that $N < q < \infty$ in view of 
the Sobolev imbedding theorem: $\|v\|_{L_\infty(\Omega)}
\leq C\|v\|_{W^1_q(\Omega)}$ with $\Omega=\BR^N_\pm$ and 
$\Omega=\BR^N$. 
To solve problem \eqref{stokeseq:1}, 
we use the maximal $L_p$-$L_q$ regualrity for the parabolic
equations.  From this point of view, we represent 
$\hrho_+$ by the integration along the characteristic
curve generated by $\bv_+$\footnote[3]{Tani \cite{T1} represented the mass
density with the help of the velocity field to prove the 
local well-posedness of the Navier-Stokes equations describing
the compressible viscous fluid flow (cf. also \cite{ST, T2}) 
It was also suggested by J. Pr\"uss to the author to represnt 
$\hrho_+$ by $\hbu_+$ and $H$ with the help of the equation of
balance of mass when the author visited Halle university in 
the early of April, 2014.}
to eliminate $\hrho_+$ from the first equation of \eqref{stokeseq:1},
which is the hyperbolic equation for $\hrho_+$.

Given function $f$ defined on $\BR^N_+$, the Lions
extension ${\rm Ext}[f]$ of $f$ is defined by 
$${\rm Ext}[f](x, t) = \begin{cases} f(x, t) \quad&\text{for
$x_N > 0$}, \\
3f(x', -x_N, t) - 2f(x', -2x_N, t) 
\quad&\text{for $x_N < 0$},
\end{cases}
$$
Set $\hat\bw_+ = {\rm Ext}[\hbu_+]$ and 
$\bv = (\hat w_{+1}, \ldots, \hat w_{+ N-1}, 
\hat w_{+N} - K_0 - \sum_{j=1}^N K_j\hat w_{+j})$. 
Note that  $\bv = \bv_+$ on $\BR^N_+$. 
We assume that  
\begin{equation}\label{small:2}
\int^T_0 \|\nabla\bv(\cdot, t)\|_{L_\infty(\BR^N)}\,dt
\leq \epsilon_1
\end{equation}
with some small positive constant $\epsilon_1 > 0$.
We use the usual fixed point argument to solve the nonlinear 
problem and in this argument we keep the situation
where $\hbu_+$ and $H$ satisfy \eqref{small:2}.

Let $\hat \xi$ be the solution to the Cauchy problem:
$$\frac{d}{dt}\hat \xi(\eta, t) = \bv(\hat\xi(\eta, t), t),\quad
\hat\xi(\eta, 0) = \eta \in \BR^N.
$$
According to Str\"ohmer \cite{St}, we choose an $\epsilon_1 > 0$ so 
small that the map: $\eta \mapsto \xi$ is bijective on 
$\BR^N$ for any $t \in [0, T]$.  We denote its inverse map 
by $\hat\eta = \hat\eta(\xi, t)$. 
We look for $\hrho_+$ satisfying the equation:
\begin{equation}\label{mass:1}
\pd_t\hrho_+ + \bv\cdot\nabla\hrho_+ + 
\hrho_+(\dv\hbw_+ + V_\dv(\hbw_+, H)) = 0 \quad
\text{in $\BR^N\times(0, T)$}.
\end{equation}
Since
$$\frac{\pd}{\pd t}\hrho_+(\hat\xi(\eta, t), t)
= (\pd_t\hrho_+ + \bv\cdot\nabla \hrho_+)(\hat\xi(\eta, t),t)
= g(\hat\xi(\eta, t), t)\rho_+(\hat\xi(\eta, t), t),
$$
with $g=- (\dv\hbw_+ + V_\dv(\hbw_+, H))$, 
defining $\hrho_+(\xi, t)$ by 
\begin{equation}\label{formula:mass}
\hrho_+(\xi, t)
= (\rho_{*+} + \tilde\rho_{0+}(\eta))
e^{-\int^t_0(\dv\hat\bw_++V_\dv(\hat\bw_+, H))(\hat\xi(\eta, s),s)\,ds}
\end{equation}
with $\eta = \hat\eta(\xi, t)$, where $\tilde\rho_{0+}(\eta)
={\rm Ext}
[\hrho_{0+}]$ to $\BR^N$, we see that $\hrho_+$ is a required function 
satisfying \eqref{mass:1} with $\hrho_+(\xi, 0)
= \rho_{*+} + \hrho_{0+}(\xi)$ in $\BR^N_+$. 

Inserting the formula of $\hrho_+$ given in \eqref{formula:mass}
into the right-hand sides: $\bF_+ = \bF_+(\hrho_+,
\hbu_\pm, H)$, $F_{\theta+} = F_{\theta+}(\hrho_+, \hbu_\pm,
\htheta_+, H)$, $G_j= G_j(\hrho_+,\bu_\pm, H)$
($j=1, \ldots, N+1$) and $G_\theta
=G_\theta(\hrho_+, \hbu_\pm, \htheta_\pm, H)$ in
\eqref{stokeseq:1}, we have the interface problem of the 
final form, which is a quasilinear parabolic
equation. As the linearized problem,
we have the decoupled two systems.  One is 
the Stokes equation with interface condition: 
\begin{alignat}2
&\quad \,\rho_{*+}\pd_t\bu_+ - \DV\bS_{*+}(\bu_+) = \bff_+
&\quad&\text{in $\BR^N_+\times(0, T)$}\nonumber \\
&
\left\{\begin{aligned}
&\rho_{*-}\pd_t\bu_- - \DV\bS_{*-}(\bu_-)+\nabla\pi_- = \bff_- \\
&\dv \bu_- = f_\dv = \dv\bff_\dv
\end{aligned}\right.
&\quad&\text{in $\BR^N_-\times(0, T)$}\label{stokes:1} 
\end{alignat}
subject to the interface condition: for $x \in \BR^N_0$ and $t \in (0, T)$
\begin{equation}\label{stokes:2}\begin{split}
&\mu_{*-}D_{iN}(\bu_-)|_- - \mu_{*+}D_{iN}(\bu_+)|_+ = g_i
\quad(i = 1, \ldots, N-1), \\
&(\mu_{*-}D_{NN}(\bu_-) - \pi_-)|_- = \sigma_-\Delta'H +g_N \\
&(\mu_{*+}D_{NN}(\bu_+)
+ (\lambda_{*+} - \mu_{*+})\dv\bu_+)|_+ = \sigma_+\Delta'H + g_{N+1}, \\
&u_{-i}|_- - u_{+i}|_+ = h_i \quad(i = 1, \ldots, N-1), \\
&\pd_tH - \Bigl(\frac{\rho_{*-}}{\rho_{*-}-\rho_{*+}}u_{-N}
-\frac{\rho_{*+}}{\rho_{*-} - \rho_{*+}}u_{+N}
\Bigr) = d
\end{split}\end{equation}
and the initial condition:
\begin{equation}\label{stokes:3}
\bu_\pm|_{t=0} = \bu_{0\pm} \quad\text{in $\BR^N_\pm$},\quad
H|_{t=0} = H_0 \quad\text{in $\BR^N_0$}, 
\end{equation}
where 
we have set $\sigma_\pm = \rho_\pm \sigma(\rho_{*-}-\rho_{*+})^{-1}$
and we have used the equivalent relations \eqref{eqiv:1}. 
Another is the heat equations with interface condition: 
\begin{equation}\label{heat:1}\begin{split}
\rho_{*+}\kappa_{*+}\pd_t\theta_+ - d_{*+}\Delta\theta_+
= \tilde f_+ 
&\quad\text{in $\BR^N_+\times(0, T)$} \\
\rho_{*-}\kappa_{*-}\pd_t\theta_- - d_{*-}\Delta\theta_-
= \tilde f_- 
&\quad\text{in $\BR^N_-\times(0, T)$}
\end{split}\end{equation}
subject to the interface condition: for $x \in \BR^N_0$ and $t \in (0, T)$
\begin{equation}\label{heat:2}
\theta_-|_- - \theta_+|_+=0, \quad
d_{*+}\pd_N\theta_-|_- - d_{*+}\pd_N\theta_+|_+ = \tilde g
\end{equation}
and the initial condition:
\begin{equation}\label{heat:3}
\theta_\pm|_{t=0} = \theta_{0\pm}
\quad\text{on $\BR^N_\pm$}.
\end{equation}
We have the following theorem about the maximal 
$L_p$-$L_q$ regularity for problem \eqref{stokes:1},
\eqref{stokes:2}, \eqref{stokes:3}. 
\begin{theorem}\label{thm:stokes}
Let $1 < p, q < \infty$ and $0 < T < \infty$. 
Assume that $\rho_{*-} \not= \rho_{*+}$. 
Then, for any initial data $\bu_{0\pm} \in B^{2(1-1/p)}_{q,p}(\BR^N_\pm)$ 
and $H_0 \in B^{3-1/p}_{q,p}(\BR^N)$,
and  right-hand sides
of \eqref{stokes:1} and \eqref{stokes:2} 
\begin{align*}
&\bff_\pm \in L_p((0, T), L_q(\BR^N_\pm)), \, 
f_\dv \in L_p((0, T), W^1_q(\BR^N_-)), \, \bff_\dv
\in W^1_p((0, T), L_q(\BR^N_-))\\
&d \in L_p((0, T), W^2_q(\BR^N)),\,
 g_i \in L_p((0, T), W^1_q(\BR^N)) \cap W^1_p((0, T), W^{-1}_q(\BR^N)) \\
&h_j \in L_p((0, T), W^2_q(\BR^N)) \cap W^1_p((0, T), L_p(\BR^N)) 
\end{align*}
for $i=1, \ldots, N+1$ and $j=1, \ldots, N-1$,   
satisfying the compatibility conditions:
\begin{alignat*}2
&\dv \bu_{0-} = f_-|_{t=0} = \dv\bff_\dv|_{t=0} 
&\quad&\text{in $\BR^N_-$}, \\
&\mu_{*-}D_{iN}(\bu_{0-})|_- - \mu_{*+}D_{iN}(\bu_{0+})|_+ 
= g_i|_{t=0}
\quad(i = 1, \ldots, N-1) &\quad
&\text{on $\BR^N_0$}, \\
&
(\mu_{*+}D_{NN}(\hbu_{0+}) + (\lambda_{*+}-\mu_{*+})
\dv\hbu_{0+})|_+\\
&\quad = \sigma_+\Delta'H_0 + g_{N+1}|_{t=0}
&\quad&\text{on $\BR^N_0$}, \\
&u_{0-i}|_- - u_{0+i}|_+ = h_i|_{t=0} \quad(i = 1, \ldots, N-1)
&\quad
&\text{on $\BR^N_0$}.
\end{alignat*}
then, problem \eqref{stokes:1}, \eqref{stokes:2}, \eqref{stokes:3} 
admits unique solutions $\bu_\pm$ and $H$ with
\begin{align*}
\bu_\pm & \in L_p((0, T), W^2_q(\BR^N_\pm))\cap
W^1_p((0, T), L_q(\BR^N_\pm)), \\
H & \in L_p((0, T), W^3_q(\BR^N)) \cap W^1_p((0, T), 
W^2_q(\BR^N))
\end{align*}
possessing the estimates:
\begin{multline*}
\sum_{\ell=\pm} \{\|\bu_\ell\|_{L_p((0, t), W^2_q(\BR^N_\ell))}
+ \|\pd_t\bu_\ell\|_{L_p((0, t), L_q(\BR^N_\ell))}\}\\
+ \|\pd_tH\|_{L_p((0, t), W^2_q(\BR^N))}
+ \|H\|_{L_p((0, t), W^3_q(\BR^N))} \\
\leq Ce^{\gamma t}\{\sum_{\ell=\pm}
(\|\bu_{0\ell}\|_{B^{2(1-1/p)}_{q,p}(\BR^N_\ell)}
+ \|\bff_\ell\|_{L_p((0, t), L_q(\BR^N_\ell))})
+ \|f_\dv\|_{L_p((0, t), W^1_q(\BR^N_-)}\\
+\|\bff_\dv\|_{L_p((0, T), L_q(\BR^N_-))} 
+ \sum_{i=1}^{N+1}(\|g_i\|_{L_p((0, t), W^1_q(\BR^N))}
+ \|\pd_tg_i\|_{L_p((0, t), W^{-1}_q(\BR^N))})\\
+ \sum_{j=1}^{N-1} (\|h_j\|_{L_p((0, t), W^2_q(\BR^N))}
+ \|\pd_th_j\|_{L_p((0, t), L_q(\BR^N))})
+ \|d\|_{L_p((0, t), W^2_q(\BR^N))}\}
\end{multline*}
for any $t \in (0, T)$ with some positive constants
$C$ and $\gamma$ independent of $t$ and $T$.
\end{theorem}
And also, we have the following theorem about
the maximal $L_p$-$L_q$ regularity for
 problem \eqref{heat:1},
\eqref{heat:2}, \eqref{heat:3}.
\begin{theorem}\label{thm:heat}
Let $1 < p, q < \infty$ and $0 < T < \infty$. Then, 
for any initial data $\theta_{0\pm} \in B^{2(1-1/p)}_{q,p}(\BR^N_\pm)$
and  right-hand sides
$$\tilde f_\pm \in L_p((0, T), L_q(\BR^N_\pm)),\quad
\tilde g \in L_p((0, T), W^1_q(\BR^N)) \cap W^1_p((0, T), W^{-1}_q
(\BR^N))$$
satisfying the compatibility condition:
$$[[\theta_0]]=0, \quad d_{*-}\pd_N\theta_{0-}|_- - d_{*-}\pd_N\theta_{0+}|_+
= \tilde g|_{t=0}\quad\text{on $\BR^N_0$},$$
problem \eqref{heat:1} and \eqref{heat:2} admits 
unique solutions $\theta_\pm$ with
$$\theta_\pm \in L_p((0, T), W^2_q(\BR^N_\pm)) \cap 
W^1_p((0, T), L_q(\BR^N_\pm))$$
satisfying the estimate:
\begin{multline*}
\sum_{\ell=\pm}\{\|\theta_\ell\|_{L_p((0, t), W^2_q(\BR^N_\ell))}
+ \|\pd_t\theta_\ell\|_{L_p((0, t), L_q(\BR^N_\ell))}\} \\
 \leq C^{\gamma t}\{\sum_{\ell=\pm}
(\|\theta_{0\ell}\|_{B^{2(1-1/p)}_{q,p}(\BR^N_\ell)}
+ \|\tilde f_\ell\|_{L_p((0, t), L_q(\BR^N_\ell))}) \\
+ \|\tilde g\|_{L_p((0, t), W^1_q(\BR^N))}
+ \|\pd_t\tilde g\|_{L_p((0, t), W^{-1}_q(\BR^N))})\}
\end{multline*}
for any $t \in (0, T)$ with some positive constants
$C$ and $\gamma$ independent of $t$ and $T$.
\end{theorem}
\begin{remark}
\thetag1~The proof of Theorem \ref{thm:stokes} is given in \cite{S4}.   
The proof of Theorem \ref{thm:heat} is found in \cite{EPS}, but it can be 
proved by using the same argument as in the proof of Theorem \ref{thm:stokes} 
in \cite{S4}.
\\
\thetag2~ Theorem \ref{main:2} is proved with the help of Theorem
\ref{thm:stokes} and Theorem \ref{thm:heat}, the Banach fixed 
point argument and, some bootstrap arguments.
The argument is quite standard, so that we may omit the proof of
Theorem \ref{main:2} (cf. Pr\"uss \cite{Pru03}).
\end{remark}
\section{$\CR$-bounded solution operators}
To prove Theorem \ref{thm:stokes}, we consider the following generalized 
resolvent problem:
\begin{alignat}2
&\rho_{*+}\lambda\bu_+ - \DV\bS_{*+}(\bu_+) = \bff_+
&\quad&\text{in $\BR^N_+$}\nonumber \\
&\rho_{*-}\lambda\bu_- - \DV\bS_{*-}(\bu_-)+\nabla\pi_- = \bff_-,
\quad\dv \bu_- = f_\dv = \dv\bff_\dv
&\quad&\text{in $\BR^N_-$}\label{rstokes:1} 
\end{alignat}
subject to the interface condition: for $x \in \BR^N_0$
\begin{equation}\label{rstokes:2}\begin{split}
&\mu_{*-}D_{iN}(\bu_-)|_- - \mu_{*+}D_{iN}(\bu_+)|_+ = g_i
\quad(i = 1, \ldots, N-1), \\
&(\mu_{*-}D_{NN}(\bu_-) - \pi_-)|_- =\sigma_-\Delta'H + g_N, \\
&(\mu_{*+}D_{NN}(\bu_+)
+ (\lambda_{*+} - \mu_{*+})\dv\bu_+)|_+ = \sigma_+\Delta'H + g_{N+1}, \\
&u_{-i}|_- - u_{+i}|_+ = h_i \quad(i = 1, \ldots, N-1), \\
&\lambda H - \Bigl(\frac{\rho_{*-}}{\rho_{*-}-\rho_{*+}}u_{-N}
-\frac{\rho_{*+}}{\rho_{*-} - \rho_{*+}}u_{+N}
\Bigr) = d,
\end{split}\end{equation}
which is corresponding to the time dependent problem \eqref{stokes:1},
\eqref{stokes:2}, \eqref{stokes:3}. 

Before stating the main result of this section, we first introduce the definition of
$\CR$-boundedness and the operator valued Fourier multiplier
theorem due to Weis \cite{Weis}.
\begin{definition}\label{def:2.1}
Let $X$ and $Y$ be two Banach spaces with norms
$\|\cdot\|_X$ and $\|\cdot\|_Y$, respectively. 
A family of operators $\CT \subset 
\CL(X, Y)$ is called  $\CR$-bounded on $\CL(X, Y)$, 
if there exist constants $C > 0$ and $p \in [1, \infty)$ such 
that for any $n\in \BN$, $\{T_j\}_{j=1}^n 
\subset \CT$, $\{f_j\}_{j=1}^n \subset X$ and sequences 
$\{r_j(u)\}_{j=1}^n$ of independent, symmetric, 
$\{-1, 1\}$-valued random variables on $[0, 1]$ there holds the
inequality:
$$\Bigl\{\int^1_0\|\sum_{j=1}^nr_j(u)T_jf_j\|_Y^p\,du\Bigr\}^{\frac1p}
\leq C\Bigl\{\int^1_0\|\sum_{j=1}^n r_j(u)f_j\|_X^p\,du\Bigr\}^{\frac1p}.$$
The smallest such $C$ is called $\CR$-bound of $\CT$, 
which is denoted by $\CR_{\CL(X, Y)}(\CT)$. Here, 
$\CL(X, Y)$ denotes the set of all bounded linear operators from 
$X$ into $Y$. 
\end{definition}
Let $\CD(\BR, X)$ and $\CS(\BR, X)$ be the set of all 
$X$ valued $C^\infty$ functions having compact supports and the Schwartz
space of rapidly decreasing $X$ valued functions, respectively, 
while $\CS'(\BR, X) = \CL(\CS(\BR, \BC), X)$. 
Given $M \in L_{1, {\rm loc}}(\BR\setminus\{0\},
X)$, we define the operator $T_M: \CF^{-1}\CD(\BR, X) \rightarrow
\CS'(\BR, Y)$ by 
\begin{equation}\label{ovfm}
T_M\phi = \CF^{-1}[M\CF[\phi]], \quad(\CF[\phi] \in \CD(\BR, X)),
\end{equation}
The following theorem
is obtained by Weis \cite{Weis}.
\begin{theorem}\label{Weis}
Let $X$ and $Y$ be two UMD Banach spaces and $1 < p < \infty$.
Let $M$ be a function in $C^1(\BR\setminus\{0\}, \CL(X, Y))$ such that
$$\CR_{\CL(X, Y)}(\{(\tau\frac{d}{d\tau})^\ell M(\tau) \mid
\tau \in \BR\setminus\{0\}) \leq \kappa < \infty
\quad(\ell = 0, 1)$$
with some constant $\kappa$.  Then, the operator $T_M$
defined in \eqref{ovfm} is extended to a bounded linear operator 
from $L_p(\BR, X)$ into $L_p(\BR, Y)$.  Moreover, denoting this 
extension by $T_M$, we have 
$$\|T_M\|_{\CL(L_p(\BR, X), L_p(\BR, Y))} \leq C\kappa$$
for some positive constant $C$ depending on $p$, $X$ and $Y$. 
\end{theorem}
\begin{remark}
For the definition of UMD space, we refer to a book due to 
Amann \cite{Amann2}.  For $1 < q < \infty$, Lebesgue space
$L_q(\Omega)$ and Sobolev space $W^m_q(\Omega)$ are both 
UMD spaces.  
\end{remark}
\begin{theorem}\label{main:r-bound} Let $1 < q < \infty$ and $0 < \epsilon
< \pi/2$.  Set $\Sigma_{\epsilon, \lambda_0} = \{\lambda \in \Sigma_\epsilon \mid
|\lambda| \geq \lambda_0\}$ $(\lambda_0 > 0)$ with 
$\Sigma_\epsilon = \{\lambda = \gamma + i\tau \in \BC\setminus\{0\} \mid 
|\arg\lambda| \leq \pi-\epsilon\}$, and 
\begin{align*}
&X_q = \{(\bff_+, \bff_-, f_\dv, \bff_\dv, \bg, \bh, d)\mid
\bff_+ \in L_q(\BR^N_+), \, \bff_-, \bff_\dv \in L_q(\BR^N_-), \, 
d  \in W^2_q(\BR^N), \\
& f_\dv \in W^1_q(\BR^N), \,
\bg=(g_1, \ldots, g_{N+1}) \in 
W^1_q(\BR^N), \, \bh=(h_1, \ldots, h_{N-1})\in W^2_q(\BR^N)\}, \\
&\CX_q = \{\bF=(\bF_{+1}, \bF_{-1}, F_{-2}, \bF_{-3}, \bF_{-4}, 
\bF_1, \bF_2, \bF_3, \bF_4, \bF_5, F_6) 
\mid \bF_{\pm 1} \in L_q(\BR^N_\pm), \enskip \\
&\quad F_{-2}, \bF_{-3}, \bF_{-4} \in L_q(\BR^N_-), \enskip 
\bF_1, \bF_2, \bF_3, \bF_4, \bF_5 \in L_q(\BR^N),\enskip F_6 \in W^2_q(\BR^N)\}.
\end{align*}
Then, there exist a constant $\lambda_0 > 0$ and operator families 
\begin{align*}
\CA_\pm(\lambda)&\in \Hol(\Sigma_{\epsilon, \lambda_0}, \CL(X_q, W^2_q(\BR^N_\pm))),
\enskip
\CP_- \in \Hol(\Sigma_{\epsilon, \lambda_0}, \CL(X_q, \hat W^1_q(\BR^N_-))),\\
\CH(\lambda) &\in  \Hol(\Sigma_{\epsilon, \lambda_0}, \CL(X_q, W^3_q(\BR^N)))
\end{align*}
such that
$\bu_\pm = \CA_\pm(\lambda)\bF_\lambda$, $\pi_- = \CP_-(\lambda)\bF_\lambda$
and $H = \CH(\lambda)\bF_\lambda$ 
are unique solutions of problem \eqref{rstokes:1}
and \eqref{rstokes:2}  for any $\lambda \in \Sigma_{\epsilon, \lambda_0}$ and 
$\bF = (\bff_+, \bff_-, f_\dv, \bff_\dv, \bg, \bh, d) \in X_q$, and we have 
\begin{align*}
\CR_{\CL(\CX_q, L_q(\BR^N_\pm))}(\{(\tau\pd_\tau)^\ell
G^1_\lambda\CA_\pm(\lambda) \mid \lambda \in \Sigma_{\epsilon, \lambda_0}\})
\leq c \quad(\ell = 0, 1),\\
\CR_{\CL(\CX_q, L_q(\BR^N_\pm))}(\{(\tau\pd_\tau)^\ell
\nabla\CP_-(\lambda) \mid \lambda \in \Sigma_{\epsilon, \lambda_0}\})
\leq c \quad(\ell = 0, 1), \\
\CR_{\CL(\CX_q, W^2_q(\BR^N_\pm))}(\{(\tau\pd_\tau)^\ell
G^2_\lambda\CH(\lambda) \mid \lambda \in \Sigma_{\epsilon, \lambda_0}\})
\leq c \quad(\ell = 0, 1)
\end{align*}
with some constant $c$.  Here, 
\begin{align*}
&G^1_\lambda\CA_\pm(\lambda) = 
(\lambda\CA_\pm(\lambda), \lambda^{1/2}\nabla\CA_\pm(\lambda), \nabla^2\CA_\pm
(\lambda)), \quad
G^2_\lambda \CH(\lambda) = (\lambda\CH(\lambda), \nabla\CH(\lambda)), \\
& 
\bF_\lambda = (\bff_+, \bff_-, \lambda^{1/2}f_\dv, \nabla f_\dv, 
\lambda\bff_\dv, \lambda^{1/2}\bg, \nabla\bg, \lambda\bh, 
\lambda^{1/2}\nabla\bh, \nabla^2\bh, d), \\
&\hat W^1_q(\BR^N) = \{\pi_- \in L_{q, {\rm loc}}(\BR^N) \mid 
\nabla \pi_- \in L_q(\BR^N)\},
\end{align*}
$\Hol(U, X)$ denotes the set of all
holomorphic functions defined on $U$ with their values in $X$,
$\nabla = (\pd_1, \ldots, \pd_N)$ and $\nabla^2 = (\pd_i\pd_j
\mid i, j=1, \ldots, N)$. 
\end{theorem}
\begin{remark}
\thetag1~ $\bF_{\pm1}$, $F_{-2}$, $\bF_{-3}$, $\bF_{-4}$, 
$\bF_1$, $\bF_2$, $\bF_3$, $\bF_4$, $\bF_5$ and 
$F_6$ are corresponding variables to 
$\bff_\pm$, $\lambda^{1/2}f_\dv$, $\nabla f_\dv$, $\bff_\dv$,
$\lambda^{1/2}\bg$, $\nabla\bg$, $\lambda \bh$, $\lambda^{1/2}\nabla\bh$,
$\nabla^2\bh$ and $d$, respecitively.\\
\thetag2~ The proof of Theorem \ref{main:r-bound} is given in \cite{S4}. 
\end{remark}


\begin{thebibliography}{9}  
\bibitem{Amann2}
\newblock H.~Amann, 
\newblock \emph{ Linear and Quasilinear Parabolic Problems, Vol. I},
\newblock Birkh\"auser, Basel, 1995.
\bibitem{Denisova} \newblock
 I.~V.~Denisova, 
\newblock Evolution of compressible and incompressible fluids separated
by a closed interface, 
\newblock \emph{Interface Free Bound} \textbf{2} (3) (2000), 283--313.
\bibitem{EPS} \newblock
 J.~Escher, J.~Pr\"uss, and G.~Simonett, 
\newblock Analytic solutions for a Stefan problem
with Gibbs-Thomson correction, 
\newblock \emph{J. reine angew. Math.}, \textbf{563} (1) (2003), 1--52.
\bibitem{KSS} \newblock
T.~Kubo, Y.~Shibata, and K.~Soga, 
\newblock On the $\CR$-boundedness for the two phase problem: 
compressible-incompressible model problem, 
\newblock to appear in \emph{Boundary Value Problem}
\bibitem{Pru03} 
\newblock J.~Pr\"uss, 
\newblock Maximal regularity for evolution equations in $L_p$-spaces, 
\newblock \emph{Conf. Semin. Mat. Univ. Bari}, \textbf{285}@
(2003), 1--39.
\bibitem{PSSS} 
\newblock J.~Pr\"uss, Y.~Shibata, S.~Shimizu, and G.~Simonett, 
\newblock On well-posedness of incompressible two-phase flows with 
phase transitions: the case of equal densities, 
\newblock \emph{
Evolution Equations and Control Theory} \textbf{1} (2012), 917--941. 
\bibitem{PS}
\newblock J.~Pr\"uss and S.~Shimizu, 
\newblock On well-posedness of incmpressible two-phase flows with phase
transitions: the case of non-equal densities,
\newblock \emph{J. Evol. Equ.} \textbf{12} (2012), 917--941.
\bibitem{PSW}
\newblock J.~Pr\"uss, S.~Shimizu, and M.~Wilke, 
\newblock Qualitative behaviour of incompressible two-phase flows with phase
transitions: the case of non-equal densities,
\newblock \emph{Comm. Partial Differential Equations}
\textbf{39} (7) (2014), 1236--1283.
\bibitem{Ross}
\newblock
 D.~Rossinelli, B.~Hejazialhosseini,
P.~Hadjidoukas, C.~Bekas, A.~Curioni, A.~Bertsch,
S.~Futral, S.~Schmidt, N.~Adams, P.~Koumoutsakos,
\newblock Simulations of cloud cavitation collapse,
\newblock \emph{SC'13 Proc. International Conf. on High Performance
Computing, Networking, Storage and Analysis}, ACM New York, NY,
USA 2013, Doi:10.1145/2503210.2504565

\bibitem{S0} \newblock
Y.~Shibata, 
\newblock Generalized resolvent estimates of 
the Stokes equations with first order boundary 
condition in a general domain, 
\newblock \emph{J. Math. Fluid Mech.}, \textbf{15} (1) (2013), 1--40.  
\bibitem{S1} \newblock
Y.~Shibata, 
\newblock On the $\CR$-boundedness 
of solution operators
for the Stokes equations with free boundary condition,
\newblock \emph{Diff. Integr. Equ.},  \textbf{27} (2014), 313--368.
\bibitem{S2} \newblock
Y.~Shibata 
\newblock On some free boundary problem
for the Navier-Stokes equations in the maximal 
$L_p$-$L_q$ regularity class,
\newblock submitted on Feb.4.2014.
\bibitem{S3}\newblock 
Y.~Shibata, \newblock On the 2 phase problem including
the phase transition, \newblock
\emph{Abstract for the $39^{\rm th}$ Sapporo symposium on
PDE at Hokkaido University}, 2014, 
http://www.math.sci.hokudai.ac.jp/sympo/sapporo/
\bibitem{S4}
\newblock Y.~Shibata
\newblock On the $\CR$ boundedness for the 
two phase problem with phase transition: 
compressible-incompressible model problem.
\newblock preprint.
\bibitem{St} \newblock
G.~Str\"ohmer, 
\newblock About a certain class of parabolic-hyperbolic
system of differential equations, 
\newblock \emph{Analysis} \textbf{9} (1989), 1--39.
\bibitem{ST} \newblock
V.~A.~Solonnikov and A.~Tani,
\newblock Evolution free boundary problem for equations of motion 
of viscous compressible barotropic liquid, 
\newblock in \emph{The Navier-Stokes 
equations II - theory and numerical methods (Oberwolfach, 1991)},
\newblock Lecture Notes in Math., \textbf{1530}, Springer, Berlag,
(1992), 30--55 
\bibitem{T1} \newblock
A.~Tani, 
\newblock On the free boundary problem for 
compressible viscous fluid motion, 
\newblock \emph{J. Math. Kyoto Univ.},
\textbf{21} (1981), 839--859. 
\bibitem{T2} \newblock
A.~Tani, 
\newblock Two phase free boundary problem for 
compressible viscous fluid motion, 
\newblock \emph{J. Math. Kyoto Univ.}, 
\textbf{24} (1984), 243--267.
\bibitem{Weis} \newblock
L.~Weis, 
\newblock Operator-valued Fourier multiplier 
theorems and maximal $L_p$-regularity. 
\newblock \emph{Math. Ann.}, \textbf{319} (2001), 735--758.
\bibitem{Yamamoto} \newblock K.~Yamamoto,
\newblock On the collapse behavior of bubble clouds in 
cavitation jets, 
\newblock \emph{Proceeding of $10^{\rm th}$ Pacific Rim
International Conference on Water Jet technology}, 
April 2013, Jeju, Korea.

\end{thebibliography}
\end{document}